\documentclass[10pt,reqno,a4paper]{amsart}
\usepackage{graphicx} 

\usepackage{todonotes}
\usepackage{amsmath, stackrel}
\usepackage{amsfonts}
\usepackage{mathrsfs}  
\usepackage{amssymb}
\usepackage{amsthm}
\usepackage{hyperref}
\usepackage{graphicx}
\usepackage[noadjust]{cite}
\usepackage{caption}
\usepackage{bbm}
\usepackage{tikz}
\usetikzlibrary{decorations.pathreplacing}
\usepackage[T1]{fontenc}
\usepackage[utf8]{inputenc}
\usepackage{mathtools}

\usepackage[top=3.35cm, bottom=3.35cm, left=2.5cm, right=2.5cm, headsep=0.2in]{geometry}

\usepackage{fancyhdr}

\newtheorem{theorem}{Theorem}[section]
\newtheorem{remark}[theorem]{Remark}

\newtheorem{corollary}[theorem]{Corollary}

\newtheorem{lemma}[theorem]{Lemma}

\newtheorem{definition}[theorem]{Definition}

\newcommand{\beq} {\begin{equation}}
\newcommand{\eeq} {\end{equation}}

\usepackage{verbatim}
\immediate\write18{texcount -tex -sum  \jobname.tex > \jobname.wordcount.tex}

\providecommand{\keywords}[1]
{
  \small	
  \textbf{\textit{Keywords---}} #1
}

\title{STABILITY OF THE WEAK HAAGERUP PROPERTY UNDER GRAPH PRODUCTS}
\date{}
\author{Shubhabrata Das and Partha Sarathi Ghosh}
\keywords{Weak Haagerup Property, Graph Product,  Completely Bounded Norm, Herz-Schur Multipliers, Approximation properties\\ 2010 Mathematics Subject Classification: 20F65, 46L07, 46B28, 43A35}

\begin{document}

\maketitle

\begin{abstract}
In this paper we prove that: Any graph product of finitely many groups, all of them satisfying weak Haagerup property with $\Lambda_{WH}=1$, also satisfies weak Haagerup property and as a corollary of this result we obtain that the free product of weakly Haagerup groups with $\Lambda_{WH}=1$, again has weak Haagerup property with $\Lambda_{WH}=1$.
\end{abstract}

\section{Introduction}
A group is called \textit{amenable} if it admits a sequence of functions, converging point-wise to the constant function one (called \textit{approximate identities}), which are finitely supported, and positive definite. Amenable groups form a fairly large class which is stable under taking direct products, quotients and subgroups. Finite rank free groups are not amenable.

\medskip

In \cite{haagerup1978example}, Haagerup showed existence of an approximate identity $\{ \phi_{n} \}$ on the free group $\mathbb{F}_n$, such that they are positive definite, vanishing at infinity (with respect to word metric). Akemann-Walter \cite{akemann1981unbounded} and Choda \cite{choda1983group} called a group to have the \textit{Haagerup property} if it admits a similar approximate identity as in \cite{haagerup1978example}. Clearly amenable groups have Haagerup property. Fundamental group of hyperbolic surfaces, CAT(0)-cubulated groups, groups acting properly on a tree and many others have this property. Cherix et. al \cite{cherix2001groups} provides a good survey on this topic.

\medskip

One other way amenability can be generalized is called \textit{weak amenability}, where $G$ is required to have approximate identities consisting of finitely supported functions which are, uniformly bounded in the `completely bounded'-norm ($B_2$-norm, see definition \ref{definition completely bounded}). Since a positive definite map $\phi$ on $G$ has bounded cb-norm, this notion generalizes amenability. There is a constant, canonically associated to a weakly amenable group $G$, obtained as the infimum of the cb-norms of these approximate identities, called the Cowling-Haagerup constant of $G$ and is denoted by $\Lambda_{CH}(G)$. Clearly, every amenable group has $\Lambda_{CH}=1$. Most groups mentioned above were shown to be weakly amenable with $\Lambda_{CH}=1$(\cite{haagerup1978example}, \cite{ozawa2008weak}, \cite{mizuta2008bozejko}, \cite{guentner2010weak}). For an example of a group with $\Lambda_{CH}>1$, we mention a uniform lattice $\Gamma$ in $Sp(1,n)$, where the $\Lambda_{CH}(\Gamma)=2n-1$. This group satisfies the Kazhdan's property (T), which can be seen as a strong negation of the Haagerup property (see \cite{bekka2008kazhdan},\cite{cowling1989completely}).

 \medskip

Clearly, groups like $Sp(n,1)$ (or, lattices in them) do not satisfy the Haagerup property. It is also observed that \textit{most} weakly amenable groups with $\Lambda_{CH}=1$ happens to be groups with the Haagerup property and the famous conjecture by Cowling in this regard stated: $G$ is weakly amenable with $\Lambda_{CH}=1$ if and only if $G$ has Haagerup property (section 1.3.1, \cite{cherix2001groups}). Ozawa and Popa constructed examples of a non-weakly amenable groups in \cite{ozawa2010class}, which was later shown to be Haagerup by Cornulier et. al \cite{cornulier2012proper}. The other direction of the conjecture remains open.

\begin{figure}[ht!]
    \centering
         \begin{tikzpicture}[scale=0.55]
            \scope
            \draw (-3.5,-3.5) rectangle (2.5,3);
            \draw (-5,-5) rectangle (5,5);
            \draw (-7,-7) rectangle (7,7);
            \endscope
            
            \scope
            \draw[black!30!white] (-0.5,-2.75) rectangle (2,2.25);
            \draw[dashed] (-1.5,-4.3) rectangle (4,4.3);
            \draw[dashed] (-2,-6) rectangle (6.5,6.5);
            \endscope
            
            \draw (-2.75,-0.5)  node [rotate=270,text=black,scale=0.75] {Haagerup groups};
            \draw (-4.25,0)  node [rotate=270,text=black,scale=0.75] {weakly Haagerup groups with $\Lambda_{WH}=1$};
            \draw (-6,0)  node [rotate=270,text=black,scale=0.75] {weakly Haagerup groups};
            \draw (0.75,-0.25)  node [rotate=90,text=black,scale=0.75] {Amenable groups};
            \draw (3.25,0)  node [rotate=90,text=black,scale=0.75] {weakly amenable groups with $\Lambda_{CH}=1$};
            \draw (5.75,0)  node [rotate=90,text=black,scale=0.75] {weakly amenable groups};
        \end{tikzpicture}
    \label{fig:1}
\end{figure}

A further generalisation called the \textit{weak Haagerup Property} interpolates between the Haagerup property and weak amenability. Introduced by Knudby \cite{knudby2014semigroups} (see definition \ref{definition WH}), a group is called weakly Haagerup if there is an approximate identity of vanishing at infinity maps, uniformly bounded in $B_2$-norm, on the group. A group $G$ with the weak Haagerup property also comes equipped with a constant $\Lambda_{WH}(G)$, defined similar to $\Lambda_{CH}$ and it is evident that $\Lambda_{WH}(G)\leq \Lambda_{CH}(G)$, for any $G$. An example of a group not having weak Haagerup property is given in \cite{haagerup2015weak}.

\medskip


Both amenability and weak amenability are preserved under direct products. In \cite{bozejko1993weakly},
Bozejko and Picardello proved that free product of amenable groups is weakly amenable. Whereas it is still an open question whether free product of any two weakly amenable groups is weakly amenable or not. Though as a special case Ricard and Xu had proved that free product of two weakly amenable group with $\Lambda_{CH}=1$ also has weak amenability \cite{ricard2006khintchine}. In general, Haagerup property of groups also proved to be stable under taking direct products, free products and amalgamated free products (over finite subgroups) \cite{cherix2001groups}. Knudby in \cite{knudby2016weak} showed that direct product of two weak Haagerup groups is again weakly Haagerup. In this paper we wish to study stability of the weak Haagerup property under certain group construction (\textit{graph products}).

\medskip

The \textit{graph products of groups} $G(\Gamma)$, defined by Green in her thesis \cite{green1990graph}, is a novel way to combine a collection $\{G_v\}_{v\in V(\Gamma)}$ of groups parameterised along the vertex set $V(\Gamma)$ of a finite graph $\Gamma$. Depending on the nature of the given (finite) graph, a graph product of groups can also be thought of as an interpolation between the free product (which corresponds to a graph without an edge) and the direct product (which corresponds to the complete graph on a given set of vertices) of groups (for details see section \ref{graph product definition}).      

\medskip

For a finite simplicial graph $\Gamma$, and a collection $\{G_v\}_{v\in V(\Gamma)}$ of weakly amenable groups with $\Lambda_{CH}(G_{v})=1$ for each vertex $v\in V(\Gamma)$, it was shown by Reckwerdt in \cite{reckwerdt2017weak}, that the graph product $G(\Gamma)$ is weakly amenable. Stability of Haagerup property under graph product was proved by the present authors in \cite{das2023stability}, (see also \cite{antolin2013haagerup}). The main result of this paper establishes stability of the weak Haagerup property under graph product. Note that a group with Haagerup property has $\Lambda_{WH}=1$, and a weakly amenable group $G$ with $\Lambda_{CH}(G)=1$ is also weakly Haagerup with $\Lambda_{WH}(G)=1$. Therefore, following Knudby's program in \cite{knudby2014semigroups}, proving a weakly Haagerup group with $\Lambda_{WH}=1$ to be Haagerup settles the remaining part of the Cowling's conjecture. The main result of this paper can be seen as a supportive evidence for the conjecture.

\begin{theorem}\label{main I}
    Suppose $\Gamma$ is a finite simplicial graph and $\{G_{v}\}_{v\in V(\Gamma)}$ is a collection of weakly Haagerup groups with $\Lambda_{WH}(G_v)=1$, for each $v$. Then the group $G(\Gamma)$ has the weak Haagerup property with $\Lambda_{WH}(G(\Gamma))=1$.
\end{theorem}

The graph products are in some sense generalization of free products of groups, so as a consequence of the above theorem we have the following as well.

\begin{corollary}\label{main II}
    Suppose $A$ and $B$ are two weakly Haagerup groups with $\Lambda_{WH}(G_v)=1$. Then the group $G=A\ast B$ has the weak Haagerup property.
\end{corollary}

Following Knudby \cite{knudby2014semigroups}, a necessary and sufficient condition for a group $G$ to have weak Haagerup property with $\Lambda_{WH}=1$, is the existence of a proper, symmetric function $\phi$ on $G$, which can be expressed as sum of two kernels on the group, i.e. for any $x,y \in G$,
\begin{equation}
    \phi(y^{-1}x)=\rho(x,y)+\tau(x,y) 
\end{equation}
where $\rho$ is proper, conditionally negative definite and $\tau$ is bounded positive definite. So, one way to prove theorem \ref{main I} is to come up with a proper function $\phi$ as above, on $G(\Gamma)$. Since on each $G_v$, we have a proper function $\phi_v$ which is a sum of $\rho_v$ and $\tau_v$ as above, it would be reasonable to ask if we can combine the individual kernels $\rho_v$'s into a $\rho$ and $\tau_v$'s to a $\tau$, in order to obtain a `$\phi=\rho +\tau$' on $G(\Gamma)$. Upto a modification the sum of the conditionally negative definite kernels on $G_v$'s define a kernel of the same type on $G(\Gamma)$. A similar approach was followed in \cite{das2023stability} to combine conditionally negative definite kernels $\rho_v$'s on Haagerup vertex groups $G_v$ to obtain a conditionally negative definite kernel on $G(\Gamma)$, proving it to have the Haagerup property. But here, we were unable to integrate the bounded positive definite kernels $\tau_v$'s to get a similar kernel $\tau$ on $G(\Gamma)$, and hence unable to obtain the desired $\phi=\rho+\tau$. We prove theorem \ref{main I} following a general line of arguments found in \cite{bozejko1993weakly,ricard2006khintchine,reckwerdt2017weak}. Somewhat similar strategies were followed also in proposition 12.3.5 of \cite{brown2008textrm} and the proof of the main result of \cite{ozawa2008weak}.

\section{Preliminary}\label{preliminary}


 Throughout the paper, let $G$ be a finitely generated group. A function $\phi:G\rightarrow \mathbb{C}$ is called \textit{positive definite} (or of positive type or PD) if for all $n\in \mathbb{N}$ it satisfies the following: for any choice of $n$ complex numbers $c_1,c_2,\cdots,c_n$; and $n$ elements $g_1,g_2,\cdots,g_n$ from the group $G$, $$\sum_{i,j=1}^n c_i \overline{c_j} \phi(g_j^{-1}g_i)\geq 0$$ 

Any function $\phi:G\rightarrow \mathbb{C}$ induces a kernel $k_\phi:G\times G\rightarrow \mathbb{C}$, defined by $k_\phi(g,h):=\phi(h^{-1}g)$. So in other words, $\phi$ is positive definite if and only if for all $n\in \mathbb{N}$ and any subset $\{g_1,\cdots,g_n\}\subset G$, the matrix $[k_{ij}:=k_\phi(g_i,g_j)]\in \mathbb{M}_n(\mathbb{C})$ is positive definite.\\


Similarly a symmetric function $\phi:G\rightarrow\mathbb{C}$ (i.e. $\phi(x^{-1})=\phi(x), \forall x\in G$) is called \textit{conditionally negative definite} (CND) if for all $n\in \mathbb{N}$, the following is satisfied: given any subset $\{g_1,g_2,\cdots,g_n\}\subset G$ and any set $\{c_1,\cdots,c_n\}$ of $n$ complex numbers, with $\sum_{i=1}^n c_i=0$, we have $$\sum_{i=1}^n c_i\overline{c_j}k_\phi(g_i,g_j)=\sum_{i=1}^n c_i\overline{c_j}\phi(g_j^{-1}g_i)\leq 0$$


\begin{remark}
    A kernel $k:G\times G\rightarrow\mathbb{C}$ of PD type (resp. CND type), induces a function $\phi_k$ of the PD type (resp. CND type) on the group if the kernel is $G$-invariant, i.e. for all $x,y,z\in G$, $$k(xy,xz)=k(y,z).$$
\end{remark}
\vspace{1em}
Alternately a $\phi:G\rightarrow\mathbb{C}$ can be characterized (\cite{de1989propriete}) as PD if there exists a unitary representation $\pi:G\rightarrow \mathbb{B}(\mathcal{H})$ on a Hilbert space $\mathcal{H}$ and a unit vector $\xi \in \mathcal{H}$ such that 
\begin{equation}\label{pd condition}
\phi(g)=\langle\pi(g)\xi,\xi\rangle,    
\end{equation}
another way of saying it is that the existence of a map $\alpha :G \rightarrow \mathcal{H}$ such that for any $g,h \in G$
\begin{equation}\label{pd condition II}
    \phi(h^{-1}g)=\langle \alpha(g),\alpha(h) \rangle,
\end{equation}
and CND if there exists a map $R:G\rightarrow \mathcal{H}$ into a Hilbert space $\mathcal{H}$ such that 
\begin{equation}\label{cnd condition}
  k_\phi(g,h)=\phi(h^{-1}g)=||R(g)-R(h)||^2. 
\end{equation}

\medskip

The following relationship between CND maps and PD maps is due to Schoenberg.
\begin{lemma}\label{schoenberg}\cite{schoenberg1938metric}(\cite{brown2008textrm},Theorem D.11)
    A function $\phi:G\rightarrow\mathbb{C}$ is conditionally negative definite if and only if $e^{-t\phi}$ is positive definite for every $t>0$.
\end{lemma} 

\subsection{Completely Bounded Maps and $B_2$-norm}
We start with the definition of the completely bounded norm or the $B_2$-norm of a function $\phi:G\rightarrow\mathbb{C}$, which will be used to define the weak Haagerup property, in what follows. 

\vspace{0.5cm}

Consider the Hilbert space $\ell^2(G)$, canonically associated to $G$, be the space of all square summable functions on $G$. $\mathbb{B}(\ell^2(G))$ is the space of all bounded operators on $\ell^2(G)$. Any element $A\in \mathbb{B}(\ell^2(G))$ can be considered in the matrix form $[A_{g,h}]_{g,h\in G}$, where $A_{g,h}= \langle A\delta_h,\delta_g\rangle$. A kernel $k:G\times G \rightarrow \mathbb{C}$, induces a multiplier $m_k:\mathbb{B}(\ell^2(G))\rightarrow \mathbb{B}(\ell^2(G))$ defined by $$[A_{g,h}]_{g,h\in G}\mapsto [k(g,h)A_{g,h}]_{g,h\in G},$$ called a \textit{Schur multiplier}. Given a function $\phi:G\to \mathbb{C}$, the Schur multiplier $m_{k_\phi}$ for the induced kernel $k_{\phi}$, is called the \textit{ Herz-Schur multiplier} if it is a bounded operator. Further for each $n\in \mathbb{N}$, we can define $M_{n,\phi}:\mathbb{M}_n(\mathbb{B}(\ell^2(G)))\rightarrow \mathbb{M}_n(\mathbb{B}(\ell^2(G)))$ by considering the map $$\quad [a_{ij}]\mapsto [m_{k_\phi}(a_{ij})].$$ Here $\mathbb{M}_n(\mathbb{B}(\ell^2(G)))$ gets its norm by identifying it with $\mathbb{B}((\ell^2(G))^n)$.  
\begin{definition}\label{definition completely bounded}
    A function $\phi:G\rightarrow \mathbb{C}$ is called \textit{completely bounded} if there is a $B>0$ such that the operator-norm of $M_{n,\phi}$, defined above is bounded by $B$, for each $n\in\mathbb{N}$. The $B_2$-norm (alternately the completely bounded norm or the cb-norm) of $\phi$ is
    $$||\phi||_{B_2}:=\sup_{n\in\mathbb{N}} ||M_{n,\phi}||$$
\end{definition}

Computing the $B_2$-norm of a function from the above definition is in general a difficult task. But we have the following equivalent characterisation of a completely bounded function, which is more useful. 
\begin{lemma}[\cite{ozawa2008weak}, Theorem 3]\label{cb equivalent}
    Let $\phi:G\rightarrow\mathbb{C}$ be a function. Then the following are equivalent:
    \begin{enumerate}
        \item $||\phi||_{B_2}\leq B$,
        \item the operator norm $||m_{k_\phi}||\leq B$,
        \item there exist a Hilbert space $\mathcal{H}$ and $\alpha,\beta \in \ell^\infty(G,\mathcal{H})$, such that $$k_\phi(g,h)= \langle \alpha(g),\beta(h)\rangle$$ with $\max\{||\alpha||_\infty,||\beta||_\infty\}\leq \sqrt{B}$
    \end{enumerate}
\end{lemma}

Suppose $\phi(g):=1$ for all $g\in G$, then it is easy to see $||\phi||_{B_2}=1$. Further if $\phi(g):=d(e,g)$ and $G$ is infinite, then it is not completely bounded as every completely bounded map should be bounded. Moreover any PD function is completely bounded and $||\phi||_{B_2} = |\phi(e)|$. This shows why amenability implies weak amenability.

\subsection{Weak Haagerup Property}
The notion of weak Haagerup property of a group was first conceived by Knudby, possibly in an attempt to answer the Cowling's conjecture \cite{knudby2014semigroups}. As mentioned earlier this notion generalizes both the Haagerup property and the weak amenability, giving a common platform to compare both the notions.  
\begin{definition}\label{definition WH}
    $G$ is \textit{weakly Haagerup} if there is a sequence of maps $\phi_n:G\rightarrow\mathbb{C}$ with 
    \begin{enumerate}
        \item $\phi_n$ vanishes at infinity, for each $n\in \mathbb{N}$,
        \item $\{\phi_n\}$ is an approximate identity i.e. $\lim_{n\rightarrow\infty}\phi_n(g)=1$ for all $g\in G$,
         \item there exists a constant $B>0$ such that $||\phi_n||_{B_2}\leq B$ for each $n\in \mathbb{N}$.
    \end{enumerate}
        The weak Haagerup constant, denoted $\Lambda_{WH}$ for $G$ is the infimum of all such $B>0$ for which there is a sequence $\{\phi_n\}_n$ satisfying the above properties.
\end{definition}

\medskip

In \cite{knudby2014semigroups}, Knudby provided an equivalent criterion for a weakly Haagerup group with $\Lambda_{WH}(G)=1$. From proposition 3.1 in \cite{knudby2014semigroups}, a group $G$ is weakly Haagerup with $\Lambda_{WH}(G)=1$ if and only if there exists a suitable proper map $\phi:G\rightarrow \mathbb{R}$ generating approximate identities $\{e^{-\frac{\phi}{n}}:G\rightarrow\mathbb{R}\}_{n\in \mathbb{N}}$, vanishing at infinity on $G$ with $||e^{-\frac{\phi}{n}}||_{B_2}\leq 1$. Moreover proposition 4.3 of \cite{knudby2014semigroups} shows that the above semi-group generator $\phi$ is given by two functions $R,S:G\rightarrow\mathcal{H}$ into some real Hilbert space $\mathcal{H}$ satisfying:

\begin{equation}\label{equivalent of weak haagerup}
    \phi(y^{-1}x)=||R(x)-R(y)||^2+||S(x)+S(y)||^2, \quad \text{for all}~x,y\in G.
\end{equation}

\medskip

 Moreover without loss of generality we may assume that $R(1_G)=0_\mathcal{H} \in \mathcal{H}$. Setting $x=y$ in equation \ref{equivalent of weak haagerup} notice that the image of $G$ under $S$ lies on the sphere of radius $\sqrt{\frac{\phi(1_G)}{4}}$ inside $\mathcal{H}$, i.e. for any $x\in G$:
\begin{equation*}
    ||S(x)||^2=\frac{\phi(1_G)}{4}
\end{equation*}

\begin{remark}
    By expanding the second term of the right side of equation \ref{equivalent of weak haagerup}, we get the above mentioned expression $\phi=\rho+\tau$: ~ for any $x,y \in G$ 
    \begin{equation*}
        \phi(y^{-1}x)=\underbrace{||R(x)-R(y)||^2+\frac{\phi(1_G)}{2}}_{\mbox{$\rho(x,y)$}}~~ + ~~ \underbrace{2\langle S(x),S(y)\rangle}_{\mbox{$\tau(x,y)$}}
    \end{equation*}
\end{remark}

\medskip

Lemma \ref{cb equivalent} (3), guarantees the existence of a Hilbert space $\widehat{\mathcal{H}}$, so that the above kernel can be represented as inner-products of two $\widehat{\mathcal{H}}$-valued maps from $G$. In the next section we briefly describe how to do that.

\begin{remark}\label{weak Haagerup data}
    We call the tuple $(G,\phi,R,S,\mathcal{H})$ a weak Haagerup data for a weakly Haagerup group $G$, with $\Lambda_{WH}(G)=1$.  
\end{remark}

\subsection{Exponential of a Hilbert Space}\label{Exponential Hilbert space}


The Schoenberg's lemma \ref{schoenberg}, says, for a CND map $\phi$ on $G$, the function $e^{-\phi}$ will be PD. A CND function $\phi$ has a form $\phi(y^{-1}x)=||R(x)-R(y)||^2$ (equation \ref{cnd condition}). Therefore referring to equation \ref{pd condition II}, there exist a map $\alpha:G\rightarrow \mathcal{K}$ into a Hilbert space $\mathcal{K}$, so that for any $x,y \in G$ one has
\begin{equation}\label{schoenberg proof}
    e^{-||R(x)-R(y)||^2}=\langle \alpha(x), \alpha (y) \rangle
\end{equation}

In the following we describe a standard way to prove this by considering the above $\mathcal{K}$ as the ``exponential'' of the real Hilbert space $\mathcal{H}$ associated to $R$ and $S$, so that the above equation \ref{schoenberg proof} holds (See Appendix D of \cite{brown2008textrm}, or \cite{dadarlat2003constructions}). 

\medskip

For a real Hilbert space $\mathcal{H}$, consider: 
$$Exp(\mathcal{H}):=\mathbb{C}\oplus (\bigoplus_{n\geq 1} \mathcal{H}^{\otimes n} )$$
Define the map $Exp:\mathcal{H}\rightarrow Exp(\mathcal{H})$ as:
$$\xi \mapsto 1 \oplus \xi \oplus (\frac{1}{\sqrt{2!}}\xi \otimes \xi)\oplus (\frac{1}{\sqrt{3!}}\xi \otimes \xi \otimes \xi)\oplus \cdots $$
Then for any two vectors $\xi,\eta \in \mathcal{H}$ we get 
\begin{equation*}
    \langle Exp(\xi),Exp(\eta)\rangle=e^{\langle\xi,\eta\rangle},
\end{equation*}
and scaling the above map point-wise we further get 
\begin{align*}
    Exp_o:\mathcal{H} &\rightarrow Exp(\mathcal{H})\\
    \xi &\mapsto Exp_o(\xi):= e^{-||\xi||^2}Exp(\sqrt{2}\xi)
\end{align*}
such that $\langle Exp_o(\xi),Exp_o(\eta) \rangle= e^{-||\xi-\eta||^2}$.

\medskip

Therefore for a group $G$, with $\Lambda_{WH}(G)=1$, one has  for $x,y\in G$ 

\begin{equation}
    \begin{array}{rl}\label{lemma cb equivalent computation}
         \langle Exp_o(\frac{R(x)}{\sqrt{n}}),Exp_o(\frac{R(y)}{\sqrt{n}}) \rangle & =  e^{-\frac{||R(x)-R(y)||^2}{n}}\\
         \\
        \langle Exp_o(\frac{S(x)}{\sqrt{n}}),Exp_o(-\frac{S(y)}{\sqrt{n}}) \rangle & = e^{-\frac{||S(x)+S(y)||^2}{n}}
    \end{array}
\end{equation}

\medskip

Hence we have the analogue of lemma \ref{cb equivalent} for the functions $e^{-\frac{\phi}{n}}$ for groups $G$ with $\Lambda_{WH}=1$. Observe that the first kernel in equation \ref{equivalent of weak haagerup} is CND. Lemma \ref{schoenberg}, tells us that the kernel $G\times G \ni (x,y)\mapsto e^{-||R(x)-R(y)||^2}$ is PD and therefore the $B_2$-norm of the kernel is 1. A simple calculation gives us $||Exp_o(\frac{S(x)}{\sqrt{n}})||, ||Exp_o(-\frac{S(y)}{\sqrt{n}})||$ equal to 1; therefore from lemma \ref{cb equivalent} and equation \ref{lemma cb equivalent computation} we can say that the other part of the kernel $(x,y)\mapsto e^{-\frac{\phi(y^{-1}x)}{n}}$ also has bounded $B_2$-norm i.e.: 
$$||(x,y)\mapsto e^{-\frac{||S(x)+S(y)||^2}{n}}||_{B_2}\leq 1$$

\subsection{Graph Product of Groups}\label{graph product definition}

 Let $\Gamma= (V(\Gamma),E(\Gamma))$ be a finite simplicial graph, which has no loops and no multiple edges. The \textit{graph product} $G(\Gamma)$ of a collection of groups $\{G_v:v\in V(\Gamma)\}$, indexed by the vertices of $\Gamma$, is a group generated by the groups $G_v$'s where two elements from $G_v$ and $G_w$ will commute if and only if there is an edge $[v,w]\in E(\Gamma)$ \cite{green1990graph}. 
 
 \medskip

 \begin{definition}
    Let $\{G_v\}_{v\in V(\Gamma)}$ be a collection of groups indexed by $V(\Gamma)$.  The graph product $G(\Gamma)$ of this collection is the quotient of the free product $*_{v\in V(\Gamma)}G_v$ by the normal subgroup generated by  $\{g_vg_wg_{v}^{-1}g_{w}^{-1}:g_v\in G_v,g_w\in G_w ~\text{and}~[v,w]\in E(\Gamma)\}$ i.e. $$G(\Gamma)=\frac{ *_{v\in V(\Gamma)}G_v} { \langle\langle \{[G_v,G_w]:[v,w]\in E(\Gamma)\}\rangle\rangle}$$

        \begin{figure}[ht!]
                \centering
                \parbox{5.5cm}{
                \begin{center}
                \begin{tikzpicture}[scale=0.33] 
                \draw [fill] (1,0) circle [radius=5pt] node[below]{$G_1$};
                \draw [fill] (5,0) circle [radius=5pt] node[below]{$G_2$};
                \draw [fill] (6,3.5) circle [radius=5pt] node[right]{$G_3$};
                \draw [fill] (0,3.5) circle [radius=5pt] node[left]{$G_4$};
                \draw [fill] (3,6) circle [radius=5pt] node[above]{$G_5$};
                \end{tikzpicture}
                \end{center}
                \caption{$\Gamma$ is disconnected\\ \centering{$G(\Gamma)=*_{i=1}^5G_i$}}
                \label{fig:2A}}
                \qquad
                \begin{minipage}{5.5cm}
                \begin{center}
                \begin{tikzpicture}[scale=0.33] 
                \draw [fill] (1,0) circle [radius=5pt] node[below]{$G_1$};
                \draw [fill] (5,0) circle [radius=5pt] node[below]{$G_2$};
                \draw [fill] (6,3.5) circle [radius=5pt] node[right]{$G_3$};gluing
                \draw [fill] (0,3.5) circle [radius=5pt] node[left]{$G_4$};
                \draw [fill] (3,6) circle [radius=5pt] node[above]{$G_5$};
                \draw (1,0) -- (5,0) -- (6,3.5) -- (0,3.5) -- (1,0) -- (3,6) -- (5,0) -- (0,3.5) -- (3,6) -- (6,3.5) -- (1,0);
                \end{tikzpicture}
                \end{center}
                \caption{$\Gamma$ is complete\\\centering{$G(\Gamma)=\oplus_{i=1}^5G_i$}}
                \label{fig:2B}
                \end{minipage}
        \end{figure} 
         
 \end{definition}

 Note that the graph product of groups over a graph $\Gamma$ with $n$ vertices and no edges, is the free product of the vertex groups, and if $\Gamma$ is a complete graph on $n$ vertices, the corresponding graph product turns out to be the direct product of the vertex groups, see the figures \ref{fig:2A}, \ref{fig:2B}. Considering these two cases as the two extremes (zero edges - all edges), graph product in general can be seen as an interpolation between the free product and direct product of the corresponding vertex groups. 

\medskip

Suppose $G= \ast_v G_v$ is free product of groups, then any $g \in G$ can be written as $g= g_1g_2\cdots g_n$, where each $g_i\in G_{v_i}\setminus\{e\}$ and no two consecutive $g_i$'s will come from the same group. This expression for $g$ is unique in free product of groups, and is called \textit{normal form} of $g$ \cite{lyndon1977combinatorial},\cite{serre2002trees}. 

\medskip

Similarly, if $g\in G(\Gamma)$ is an element in a graph product of groups, $g$ can also be expressed as $g=g_1g_2\dots g_n$, where $g_i\in G_{v_i}\setminus \{e\}$ and consecutive elements come from two different vertex groups. But the presence of relators $\{[G_v,G_w]\}_{[v,w]\in E(\Gamma)}$ in a general graph product $G(\Gamma)$, prevents a unique normal form expression for an element $g$ in $G(\Gamma)$. Instead each element has a well defined reduced form, which is essentially a normal form up to the `shuffles' facilitated by the relators. 

\medskip

 By a \textit{shuffle} on an expression $g=g_1g_2\cdots g_{i-1}g_ig_{i+1}\cdots g_n$ we mean that if $[v_i,v_{i+1}]\in E(\Gamma)$, then $g$ is also expressed as $g=g_1g_2\cdots g_{i-1}g_{i+1}g_i\cdots g_n$. A decomposition of $g$ as above is called \textit{reduced} if no shuffles, followed by multiplication of two consecutive elements (whenever possible), reduces the number of letters present in the decomposition.  
 
\medskip

Green proved in her thesis (Theorem 3.9, \cite{green1990graph}) that for each element $g$ in a graph product $G(\Gamma)$, the set of letters appearing in a reduced expression is unique. Therefore, the number of elements in each expression of $g$ is the same, and this induces a ``reduced length" metric on the graph product $G(\Gamma)$, given by $d_r(g,h):=|h^{-1}g|_r$, where $|g|_r$ denotes the number of letters appearing in a reduced expression of $g$.  Note that it is the normal length in the case of free products of groups.

\medskip

Let $v$ be a vertex of $\Gamma$. The \textit{star of} $v$, denoted by $st(v)$, is defined as the collection of vertices which are at a distance $\leq 1$ from $v$ in $\Gamma$, i.e. $st(v)=\{w\in \Gamma:[v,w]\in E(\Gamma)\}\cup \{v\}$. Let $\Gamma_{st(v)}\subseteq \Gamma$ be the maximum sub-graph in $\Gamma$ having vertex set $st(v)$. Define $G(st(v)):=G(\Gamma_{st(v)})$. 

\medskip

Notice that, if $g_1g_2\cdots g_n$ is a reduced expression and for $1\leq i<j\leq n$, $g_i,g_j$ are from same vertex group $G_v$, then there has to be some $k$ with $i<k<j$ such that $g_k\notin G(st(v))$. Otherwise it would contradict the definition of reduced expression.

\medskip

By a reduced form of an element $g\in G(\Gamma)$ we mean a class of reduced decompositions of $g$. The ``$d$-tail'' of $g$ is the collection of  strings of letters that appear in the last $d$-length portion of elements from the above class. 
\begin{lemma} (see lemma 2.5 in \cite{reckwerdt2017weak}\label{d tail lemma})
    Suppose that $g,h \in G(\Gamma)$ and $d = |h^{-1}g|_r$. Then each term of a reduced form of $h^{-1}g$ is either a term from the $d$-tail of $g$ or $h$, or an amalgamation of terms from the $d$-tails of $g$ and $h$.
\end{lemma}

Therefore $d$-tails of $g$ in general can contain more than $d$ letters (in free product case it is exactly $d$).The following lemma, which is obvious from discussions here, has been used in the estimates of the main theorem.
\begin{lemma}\cite{reckwerdt2017weak}\label{maximum terms}
    Suppose $M$ is maximum such that there is a complete sub-graph in $\Gamma$ with $M$ number of vertices. Then the maximum number of letters in $d$-tail of any element of $G(\Gamma)$ is $dM$.
\end{lemma}

\subsection{Construction of Cube Complex from Graph Products}\label{construction of wall space}

In order to prove the weak Haagerup property, our goal is to construct an approximate identity $\{\phi_n\}$ such that, 
\begin{itemize}
    \item each $\phi_n$ vanishes at infinity,
    \item and the $B_2$-norm of each $\phi_n$ is uniformly bounded.
\end{itemize}

In section \ref{weak haagerup property for graph product}, we produce such functions by gluing similar functions $\phi_{n,v}$'s coming from each of the vertex group $G_v$'s along a certain `base space' $X$. In both the cases \cite{bozejko1993weakly} and \cite{reckwerdt2017weak}, in order to guarantee the uniform boundedness of the approximate identities in the $B_2$-norm, the corresponding prescriptions required scaling by a suitable PD factor. We too will need a similar positive definite scaling in order to keep $||\phi_n||_{B_2}$'s uniformly bounded.

\medskip

Similar to \cite{reckwerdt2017weak}, the base space in our case is a CAT(0)-cube complex canonically associated to the graph product $G(\Gamma)$, and the distance function on a CAT(0)-cube complex was proven to be CND by Niblo and Reeves \cite{niblo1997groups}. We obtain a PD scaling factor using this CND distance function following lemma \ref{schoenberg}. In this subsection we will briefly discuss the CAT(0)-cube complex structure associated to a graph product $G(\Gamma)$.

\medskip
\medskip

 A \textit{cube complex} is a combinatorial object, obtained by gluing euclidean cubes of different dimensions, by isometries along faces of lower dimensions. The euclidean cubes induce a metric on the cube complex. A metric space $X$ is called `CAT(0)', if it satisfies the CAT(0)-comparison property i.e. every geodesic triangle in $X$ is as \textit{thin} as its euclidean comparison triangle. By a \textit{CAT(0)-cube complex} we mean a cube complex which satisfies the CAT(0) property. For details and other equivalent definitions see \cite{gromov1987hyperbolic},\cite{bestvina2014geometric}. 

 \medskip

 One can also endow a space $X$ with a CAT(0)-cube complex structure if it admits a certain ``\textit{wall space}'' structure, given by the (combinatorial) intersection data of the hyperplanes partitioning $X$ into two halves, called \textit{walls}. Informally speaking, a wall space is a space with a family of intersecting walls. 
   
 \begin{definition}(Wall Space)
    A wall space structure on $X$ is given by a pair $(X,\mathcal{W})$, where the collection $\mathcal{W}\subseteq \mathcal{P}(X)~\text{(the power set of $X$)}$ is considered as the space of \textit{half-spaces}, such that
\begin{enumerate}
    \item $\mathcal{W}$ is closed under taking complement i.e. $h\in \mathcal{W}\iff h^c\in \mathcal{W}$ 
    \item a \textit{wall}, denoted simply by $h$, is defined to be the pair $(h,h^c)$, which is uniquely determined by the element $h\in \mathcal{W}$  
    \item two points $x,y\in X$ are \textit{separated} by a wall $h$ if $x\in h$ and $y \in h^c$ 
    \item  two walls $u$ and $v$ are meant to be \textit{crossing} each other if all the four sets $u\cap v, u\cap v^c, u^c\cap v, u^c\cap v^c$ are non-empty.
\end{enumerate}
\end{definition}

In \cite{chatterji2005wall}, Chatterji and Niblo showed that if a group $G$ acts on a wall space $(X,\mathcal{W})$, respecting the wall structure (i.e. sending walls to walls) given by $\mathcal{W}$, then there is a dual CAT(0)-cube complex $\hat{X}_{_\mathcal{W}}$ associated to $(X,\mathcal{W})$, on which $G$ acts suitably.

\begin{theorem} \cite{chatterji2005wall}\label{walls}
    Suppose $G$ is a group acting on a wall space $X$, and $\mathcal{W}$ is the set of walls with an upper bound on the number of mutually crossing walls. Then there is a finite dimensional $CAT(0)$ cube complex $\hat{X}_{_\mathcal{W}}$, on which $G$ admits an isometric action.
\end{theorem}

\begin{remark}
    Although Chatterji and Niblo proved the above theorem in the proper setting (i.e. the group acting properly on the wall space), the theorem holds true both with or without the properness assumption.
\end{remark}

\subsubsection{Wall Space structure on Graph Products}\label{Wall Space structure on Graph Product}
 Let $v$ be a vertex in $\Gamma$. Define $W_v\subset G(\Gamma)$ to be the collection of those elements which have a reduced form that starts with elements from $G_v$ i.e. $$W_v=\{g \in G(\Gamma):g=g_1g_2\cdots g_m~\text{is a reduced form and}~g_1\in G_v\}.$$

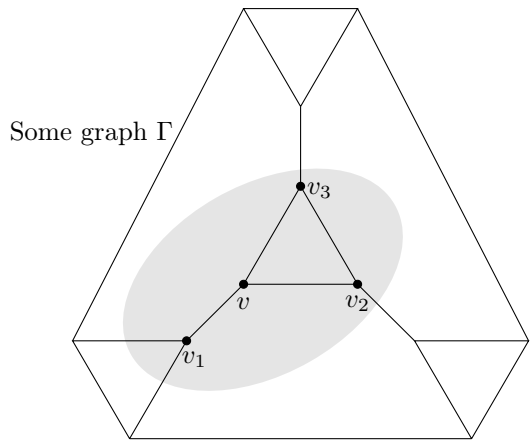
\begin{figure}[ht!]
    \centering
       \begin{tikzpicture}[scale= 0.5]
            \path [draw=none,fill=gray, fill opacity = 0.2,rotate=30] (0.5,-0.15) ellipse (4cm and 2.5cm);
            \draw  (0,0) -- (60:3) -- (3,0) -- cycle;
            \draw (0,0) -- (-1.5,-1.5);
             \draw (3,0) -- (4.5,-1.5);
             \draw (1.5,5.19/2) -- (1.5,5.19/2+2.12);
             \draw  (1.5,5.19/2+2.12) --  (3,5.19+2.12) --  (0,5.19+2.12) --(1.5,5.19/2+2.12);
             \draw  (1.5-4.5,-8.81+5.19/2+2.12) --  (-4.5+3,5.19+2.12-8.81) --  (0-4.5,5.19+2.12-8.81) --(1.5-4.5,5.19/2+2.12-8.81);
             \draw  (1.5+4.5,5.19/2+2.12-8.81) --  (3+4.5,5.19+2.12-8.81) --  (4.5,5.19+2.12-8.81) --(1.5+4.5,5.19/2+2.12-8.81);
             \draw (3,5.19+2.12)--(3+4.5,5.19+2.12-8.81);
             \draw (1.5+4.5,5.19/2+2.12-8.81)--(1.5-4.5,-8.81+5.19/2+2.12);
             \draw (0-4.5,5.19+2.12-8.81)--(0,5.19+2.12);
              \node at (0,-0.5) {$v$};
              \draw [fill] (0,0) circle [radius=3pt];
              \node at (-4,4) {Some graph $\Gamma$};
              \draw [fill] (-1.5,-1.5) circle [radius=3pt];
              \node at (-1.3,-2) {$v_1$};
             \draw [fill] (3,0) circle [radius=3pt];
              \node at (3,-0.5) {$v_2$};
            \draw [fill] (1.5,2.598) circle [radius=3pt];
              \node at (2,2.5) {$v_3$};
        \end{tikzpicture}
    \caption{$G_v,G_{v_1},G_{v_2},G_{v_3}$ generates $G(\Gamma_{st(v)})$  in $G(\Gamma)$}
    \label{fig:3}
\end{figure}

 Observe that $W_v$ is a wall in $G(\Gamma)$ (since $(W_v,W_v^c)$ forms a partition of $G(\Gamma)$). In the figure \ref{fig:3}, the star of $v$ is the sub-graph $\Gamma_{st(v)}$ of $\Gamma$, consisting of vertices $\{v,v_1,v_2,v_3\}$ along with the concerned edges and $W_v$ is exactly $G(\Gamma_{st(v)})$. Notice that $W_v^c$ is the set of elements whose reduced words start with element from the groups associated to vertices lying outside the shaded region.

 \medskip
 
 For any $g\in G(\Gamma)$, the sets $gW_v$ and $gW_v^c$ is naturally defined by left multiplication. One can see that the set $gW_v$ contains all the elements $k\in G(\Gamma)$ whose $g^{-1}$ translate has a reduced form starting with an element from $G_v$. Further it is not hard to check that $gW_v^c=(gW_v)^c$. Therefore $gW_v$ is also a wall of $G(\Gamma)$. Let $\mathcal{W}$ be the collection of walls given by 
 \begin{equation}\label{definition of walls}
    \mathcal{W}:=\{gW_v:g\in G(\Gamma), v\in \Gamma\}    
 \end{equation}

In \cite{reckwerdt2017weak}, it was showed that  $\mathcal{W}$ has an upper bound on the number of mutually crossing walls, thus one have a finite dimensional dual $CAT(0)$-cube complex $\hat{X}_{_\mathcal{W}}$, as in theorem \ref{walls}. In the following we denote $\hat{X}_{_\mathcal{W}}$ simply by $X$. The vertex set $V(X)$ is the collection of all (left)-cosets of $G(st(v))$'s in $G(\Gamma)$, for each $v\in \Gamma$, on which $G$ admits a canonical action.  Moreover there is an $x_0\in X$ such that for all $ g,h\in G(\Gamma)$ (lemma 3.5, \cite{reckwerdt2017weak}) 
\begin{equation}\label{CAT}
    d_X(gx_0,hx_0)=2|h^{-1}g|_r
\end{equation}

\medskip

The distance metric on a finite dimensional $CAT(0)$-cube complex is of CND type \cite{niblo1997groups}. The equation \ref{CAT} tells us the reduced metric on $G(\Gamma)$ is proportional to the combinatorial metric of the dual CAT(0)-cube complex, obtained from the above wall structure on $G(\Gamma)$.  Hence the reduced length function $g\mapsto |g|_r$ on $G(\Gamma)$ is a CND function.

\section{Weak Haagerup property for Graph Product}\label{weak haagerup property for graph product}

Let $\Gamma$ be a finite simplicial graph, and $\{G_v\}_{v\in V(\Gamma)}$ be a collection of groups having weak Haagerup property with weak Haagerup constant $\Lambda_{WH}(G_v)=1$. For each vertex $v\in V(\Gamma)$, we are given a weak Haagerup data $(G_v,\phi_v,R_v,S_v,\mathcal{H}_v)$. In this section we will combine this entire vertex data $\{(G_v,\phi_v,R_v,S_v,\mathcal{H}_v)\}_{v \in V(\Gamma)}$ to prove weak Haagerup property for the graph product $G(\Gamma)$ of the groups $\{G_v\}_{v\in V(\Gamma)}$. 

\medskip

From the above data, approximate identity on $G_v$ is given by $\{e^{-\frac{\phi_v}{n}}\}$. We scale $e^{-\frac{\phi_v}{n}}$ suitably to get a function $\psi_{n,v}$, in order to define 1 at the corresponding identity of the groups $G_v$, for $v \in V(\Gamma)$. Following \cite{bozejko1993weakly}, we will define an approximate identity on $G(\Gamma)$ by combining $\psi_{n,v}$'s.

\medskip

Consider $\psi_{n,v}:G_v \rightarrow \mathbb{R}$ given by $g \mapsto e^{\frac{\phi_v(1_v)}{n}}e^{-\frac{\phi_v(g)}{n}}$. Observe that $\{\psi_{n,v}\}$ is also an approximate identity of the same type as $e^{-\frac{\phi_v}{n}}$, and the $B_2$-norm: $||\psi_{n,v}||_{B_2}\leq e^{\frac{\phi_v(1_v)}{n}}$. Therefore for any $x,y \in G_v$ we have,
\begin{equation}\label{group invariance}
    \psi_{n,v}(y^{-1}x)=e^{-\frac{-\phi_v(1_v)+\phi_v(y^{-1}x)}{n}}=e^{-\frac{\big[-\phi_v(1_v)+||R_v(x)-R_v(y)||^2+||S_v(x)+S_v(y)||^2\big]}{n}}
\end{equation}

\medskip

Further from the discussion of section \ref{Exponential Hilbert space}, one can write the $\psi_{n,v}(y^{-1}x)$ as an inner-product between two vectors in a Hilbert space (see \ref{lemma cb equivalent computation}):

\begin{multline*}
    \psi_{n,v}(y^{-1}x)=\langle Exp_o(\frac{R_v(x)}{\sqrt{n}}) \otimes e^{\frac{\phi_v(1_v)}{2n}} Exp_o(\frac{S_v(x)}{\sqrt{n}}),
    \\ Exp_o(\frac{R_v(y)}{\sqrt{n}}) \otimes e^{\frac{\phi_v(1_v)}{2n}} Exp_o(-\frac{S_v(y)}{\sqrt{n}})\rangle
\end{multline*} 
where the vector norm of each component is bounded by $\sqrt{e^{\frac{\phi_v(1_v)}{n}}}$. 
For simplicity of notation we write 

\begin{equation*}
    \alpha_{n,v}(x) := \underbrace{Exp_o(\frac{R_v(x)}{\sqrt{n}})}_{\mbox{$\alpha_{n,v,R_v}(x)$}} ~\otimes~  \underbrace{e^{\frac{\phi_v(1_v)}{2n}} Exp_o(\frac{S_v(x)}{\sqrt{n}})}_{\mbox{$\alpha_{n,v,S_v}(x)$}}
\end{equation*}
and 
\begin{equation*}
    \beta_{n,v}(x) :=  \underbrace{Exp_o(\frac{R_v(x)}{\sqrt{n}})}_{\mbox{$\beta_{n,v,R_v}(x)$}} ~\otimes~  \underbrace{e^{\frac{\phi_v(1_v)}{2n}} Exp_o(\frac{-S_v(x)}{\sqrt{n}})}_{\mbox{$\beta_{n,v,S_v}(x)$}}
\end{equation*}
where $\alpha_{n,v}$ and $\beta_{n,v}$ both are maps from $G_v$ to the Hilbert space $Exp(\mathcal{H}_v)\otimes Exp(\mathcal{H}_v)$. Note that here $\alpha_{n,v,R_v}=\beta_{n,v,R_v}$.

\medskip

Let $X$ be the finite dimensional CAT(0)-cube complex underlying the graph product $G(\Gamma)$ described in the section \ref{construction of wall space}. In \cite{das2023stability}, the authors provided a way to combine $\{R_v\}$ along $X$. Here, we briefly recall that construction.

\medskip

The vertex set $V(X)$ is the collection of cosets of $G(st(v))$, for all $v\in V(\Gamma)$ i.e. $V(X)=\{G(\Gamma)/G(st(v))\}_{v\in V(\Gamma)}$. Consider the Hilbert space $\mathcal{H}_1:=\oplus_{t\in X} \mathcal{H}_t$, where $\mathcal{H}_t$ is $\mathcal{H}_v$ if $t$ is a coset $gG(st(v))$. A non-trivial $\gamma\in G(\Gamma)$ has a reduced form  $\gamma=\gamma_1\gamma_2\ldots \gamma_m$, where $\gamma_i\in G_{v_i}\leq G(st(v_i))$ for $i=1,2,\dots,m$. Define $R_\Gamma:G(\Gamma)\rightarrow \mathcal{H}_1$ by:
\begin{equation}\label{definition of R_Gamma}
    R_\Gamma(\gamma):=
        \begin{cases}
            \oplus_{i=1}^m R_{v_i}(\gamma_i)_{\gamma_1\gamma_2\ldots \gamma_{i-1}G(st(v_i))}, ~~\text{if}~\gamma\neq 1_{G(\Gamma)}\\
            0_{\mathcal{H}}~~~~~~~~~~~~~~~~~~~~~~~~~~~~~~~~~~\text{if}~\gamma=1_{G(\Gamma)}\\
        \end{cases}
\end{equation}
Then we consider the map from $G(\Gamma)$ to the exponential Hilbert space $Exp(\mathcal{H}_1)$, given by $$ \gamma \mapsto Exp_o(\frac{R_\Gamma(\gamma)}{\sqrt{n}})$$
satisfying the following property: for any $\gamma,\eta\in G(\Gamma)$

\begin{equation}\label{schoenberg in graph product}
    \langle Exp_o(\frac{R_\Gamma(\gamma)}{\sqrt{n}}),Exp_o(\frac{R_\Gamma(\eta)}{\sqrt{n}})\rangle=e^{-\frac{||R_\Gamma(\gamma)-R_\Gamma(\eta)||^2}{n}}.
\end{equation}

\medskip

Let us choose $\epsilon>0$. For each $v \in V(\Gamma)$, $G_v$ is weakly Haagerup with $\Lambda_{WH}(G_v)=1$. From the definition \ref{definition WH} and the above defined approximate identity $\psi_{n,v}$ on $G_v$, we can choose $n$ sufficiently large so that
$||\psi_{n,v}||_{B_2}\leq 1+\epsilon$. By abuse of notation, we denote that tail of the above sequence to be $\phi_{n,v}$ for which the $B_2$-norm is less than $1+\epsilon$ for each $n$. More precisely for any $x \in G_v$ we have $||\alpha_{n,v,R_v}(x)||=1=||\beta_{n,v,R_v}(x)||$. Hence
\begin{equation*}                   
||\psi_{n,v}||_{B_2}=||\alpha_{n,v,S_v}||_\infty\cdot ||\beta_{n,v,S_v}||_\infty\leq e^{\frac{\phi_v(1_v)}{n}}\leq 1+\epsilon 
\end{equation*}
and this choice of $\epsilon$ and $n$ also gives:
\begin{equation}\label{S_v estimate}
\text{for any }~ x \in G_v
    \begin{cases}
    ||\alpha_{n,v,S_v}(x)|| \leq \sqrt{1+\epsilon}\\
    ||\beta_{n,v,S_v}(x)|| \leq \sqrt{1+\epsilon}
    \end{cases}
\end{equation}

\begin{remark}\label{kernel at x,x}
    Note that for each $x\in G_v$ we have $\langle\alpha_{n,v,S_v}(x),\beta_{n,v,S_v}(x)\rangle=1$.
\end{remark}

Given the choice of $ \epsilon >0$, we consider the tuple $$\Sigma_{v,\epsilon}:=(G_v,\{\psi_{n,v}\},\{\alpha_{n,v}\},\{\beta_{n,v}\},Exp(\mathcal{H}_v)\otimes Exp(\mathcal{H}_v))$$ as an $\epsilon$-perturbed weak Haagerup data for $G_v$. In order to prove weak Haagerup  property for the graph product, we produce a similar $\epsilon$-data on $G(\Gamma)$, denoted $\Sigma_{_{\Gamma,\epsilon}}$, obtained from the entire vertex data $\{\Sigma_{v,\epsilon}\}_{v\in V(\Gamma)}$. Finally, we let the parameter $\epsilon$ tend to zero, in order to obtain $\Lambda_{WH}(G(\Gamma))=1$. 

\medskip

So far we have managed to combine the functions $\alpha_{n,v,R_v}$ and $\beta_{n,v,R_v}$ over all $v \in V(\Gamma)$ into the corresponding functions for $G(\Gamma)$ (see equation \ref{schoenberg in graph product}). To construct a desired approximate identity on $G(\Gamma)$, what remains is to glue the functions $\alpha_{n,v,S_v}$ and $\beta_{n,v,S_v}$ for all $v \in V(\Gamma)$. We break up the remainder of the proof into the following three parts.

\medskip

\begin{enumerate}
    \item[\textit{Step-I}:]  We construct a family of kernels $\psi_{n,\Gamma,d}$ (not on the entire $G(\Gamma)\times G(\Gamma)$, but) restricted to the $d$-spheres of each element of $G(\Gamma)$, with respect to the reduced length. Finally, the kernels $\psi_{n,\Gamma,d}$ will be glued together over varying $d\geq 0$, obtaining a family of kernels $\{\psi_{n,\Gamma}\}_{_{n\in \mathbb{N}}}$ on the whole group $G(\Gamma)$.
    \medskip
    \item[\textit{Step-II}:] The $B_2$-norms of the above kernels $\psi_{n,\Gamma}$ are estimated (see equation \ref{1st B_2 estimate}) by comparing their distances from suitable PD kernels $\sigma_{n,\Gamma}$ on $G(\Gamma)$ (constructed in \ref{positive deinite comparison}).
    \item[\textit{Step-III}:] We prove the kernels $\psi_{n,\Gamma}$ to be $G(\Gamma)$-invariant, and thus providing a desired approximate identity $\{\phi_{n,\Gamma}\}$ on $G(\Gamma)$.
\end{enumerate}

\medskip

\subsection{Construction of kernels}\label{section I}

First we consider an $\epsilon$-perturbed average of the vectors $\alpha_{n,v,S_v}(x)$ and $\beta_{n,v,S_v}(x)$, so that the resulting vector lies very close to the other two vectors, i.e. for any $x \in G_v$, define $$avg_{n,v,S_v}(x):=\frac{\alpha_{n,v,S_v}(x)+\beta_{n,v,S_v}(x)}{2+\sqrt{2\epsilon}}.$$ Now extend the Hilbert space $Exp(\mathcal{H}_v)$ to $\widehat{Exp(\mathcal{H}_v)}:=Exp(\mathcal{H}_v)\oplus \mathbb{C}^2\oplus \mathbb{C}^2$ and consider the map $\theta_{n,v,S_v}:G_v \rightarrow \widehat{Exp(\mathcal{H}_v)}$ given by 
\begin{align*}
    x\mapsto
        \begin{cases}
            \Bigg(avg_{n,v,S_v}(x),\begin{pmatrix}
                D_{n,v,S_v}(x)\\D_{n,v,S_v}(x)
            \end{pmatrix},\begin{pmatrix}
                0\\0
            \end{pmatrix}\Bigg), ~\text{if}~x\neq 1_v\\ 
             \\
                \Bigg(avg_{n,v,S_v}(1_v),\begin{pmatrix}
                0\\0
            \end{pmatrix},\begin{pmatrix}
                D_{n,v,S_v}(1_v)\\D_{n,v,S_v}(1_v)
            \end{pmatrix} \Bigg), ~\text{if}~ x=1_v
        \end{cases}
\end{align*}
where $D_{n,v,S_v}(x)=\sqrt{\frac{1-||avg_{n,v,S_v}(x)||^2}{2}} \in \mathbb{C}$.

\begin{remark}
    Note that from the definition of $\theta_{n,v,S_v}$ it follows that for any $x\in G_v$, $||\theta_{n,v,S_v}(x)||=1$.
\end{remark}

Let $X$ be the finite dimensional CAT(0)-cube complex on which $G(\Gamma)$ acts, (see \ref{construction of wall space}). We know that $V(X)=\cup_{v \in \Gamma}\{G(\Gamma)/G(st(v))\}$. We consider the larger Hilbert space $$\mathcal{H}_2=\otimes_{t\in V(X)}\widehat{Exp(\mathcal{H}_t)}$$  where $\mathcal{H}_t:=\mathcal{H}_v$, if $t=kG(st(v))$ for some $v$, and the corresponding vacuum vector be $\Bigg(avg_{n,v,S_v}(1_v),\begin{pmatrix}
                0\\0
            \end{pmatrix},\begin{pmatrix}
                D_{n,v,S_v}(1_v)\\D_{n,v,S_v}(1_v)
            \end{pmatrix} \Bigg)$. The reason for taking these extra $\mathbb{C}^2$ components is to absorb some ``errors'' accumulated in the process (see expressions in \ref{error manipulation} and afterwards).


\medskip 

Let $\gamma=\gamma_1\gamma_2\cdots\gamma_m$ be a reduced form for the element $\gamma\in G(\Gamma)$, and the collection of ordered $d$-tuples $(\gamma_{m-d+1},\cdots,\gamma_m)$'s appearing in the last $d$-terms of any reduced form of $\gamma$ is the $d$-tail of $\gamma$ (see section \ref{graph product definition}). For each $d\in \mathbb{N}$, we define $\alpha_{n,\Gamma,d} $ and $ \beta_{n,\Gamma,d}:G(\Gamma)\rightarrow \mathcal{H}_2$, given by

\begin{equation}\label{alpha}
\alpha_{n,\Gamma,d}(\gamma)=\otimes_{i=1}^m 
\begin{cases}
\theta_{n,\ast,S_\ast}(\gamma_i)_{\gamma_1\cdots\gamma_{i-1}G(st(v_i))}, ~~\qquad\qquad\qquad\qquad \text{if}~\gamma_i\notin d \text{-tail of}~\gamma,\\
\\
\Bigg(\alpha_{n,\ast,S_\ast}(\gamma_i),\begin{pmatrix}
             C^{\alpha_{n,\ast,S_\ast}}(\gamma_i,\gamma_i)\\0
         \end{pmatrix},\begin{pmatrix}
             C^{\alpha_{n,\ast,S_\ast}}(\gamma_i,1_\ast)\\0
         \end{pmatrix}\Bigg)_{\gamma_1\cdots\gamma_{i-1}G(st(v_i))},\\
         \qquad\qquad\qquad\qquad\qquad\qquad\qquad\qquad\qquad\qquad\text{if}~\gamma_i\in d \text{-tail of}~\gamma
\end{cases}
\end{equation}

\smallskip

\begin{equation}\label{beta}
\beta_{n,\Gamma,d}(\gamma)=\otimes_{i=1}^m
\begin{cases}
\theta_{n,\ast,S_\ast}(\gamma_i)_{\gamma_1\cdots\gamma_{i-1}G(st(v_i))}, ~~\qquad\qquad\qquad\qquad \text{if}~\gamma_i\notin d \text{-tail of}~\gamma,\\
\\
\Bigg(\beta_{n,\ast,S_\ast}(\gamma_i),
        \begin{pmatrix}
             0\\C^{\beta_{n,\ast,S_\ast}}(\gamma_i,\gamma_i)
         \end{pmatrix},
         
         \begin{pmatrix}
             0\\C^{\beta_{n,\ast,S_\ast}}(\gamma_i,1_\ast)
         \end{pmatrix}\Bigg)_{\gamma_1\cdots\gamma_{i-1}G(st(v_i))},\\
         \qquad\qquad\qquad\qquad\qquad\qquad\qquad\qquad\qquad\qquad\text{if}~\gamma_i\in d \text{-tail of}~\gamma
\end{cases}
\end{equation}

\noindent where $\ast$ stands for that $v \in \Gamma$, such that $\gamma_i \in G_v$ and the expressions for $ C^{\alpha_{n,v,S_v}}(x,y)$ and $C^{\beta_{n,v,S_v}}(x,y)$ are given by: for all $x,y \in G_v$, 

 \begin{align*}
    C^{\alpha_{n,v,S_v}}(x,y) &:=\frac{\langle\alpha_{n,v,S_v}(x),\beta_{n,v,S_v}(y)-avg_{n,v,S_v}(y)\rangle}{D_{n,v,S_v}(y)}\\
      C^{\beta_{n,v,S_v}}(x,y) &:=\frac{\langle\alpha_{n,v,S_v}(y)-avg_{n,v,S_v}(y),\beta_{n,v,S_v}(x)\rangle}{D_{n,v,S_v}(y)}
 \end{align*}

Now we prove the well-definedness of the above maps. Suppose there is an edge between $v_i$ and $v_{i+1}$ vertices of $\Gamma$. Then we have 
$$\gamma=\gamma_1\gamma_2\cdots\gamma_{i-1}\gamma_i\gamma_{i+1}\cdots\gamma_m=\gamma_1\cdots \gamma_{i-1}\gamma_{i+1}\gamma_i\gamma_{i+2}\cdots \gamma_m$$
Checking the well-definedness, then reduces to validation of the following two equations:
    \begin{align*}
        \gamma_1\cdots\gamma_{i-1}G(st(v_i))=\gamma_1\cdots\gamma_{i-1}\gamma_{i+1}G(st(v_i)) \\
        \gamma_1\cdots\gamma_{i}G(st(v_{i+1}))=\gamma_1\cdots \gamma_{i-1}G(st(v_{i+1}))
    \end{align*}
This is true due to the fact that $[G_{v_i},G_{v_{i+1}}]=1$ in $G(\Gamma)$. Since the expression of the reduced form for $\gamma$ is unique up to a finite number of shuffles between two consecutive factors whenever they commute. So the maps $\alpha_{n,\Gamma,d}$ and $\beta_{n,\Gamma,d}$ are well-defined. 

\medskip

Therefore we can define a kernel 
\begin{equation*}
    \begin{array}{rl}
       \psi_{n,\Gamma,d}:G(\Gamma)\times G(\Gamma)  &\rightarrow \mathbb{C}  \\
       (\gamma,\eta)  & \mapsto \langle\alpha_{n,\Gamma,d}(\gamma),\beta_{n,\Gamma,d}(\eta)\rangle
    \end{array}
\end{equation*}
Now set $E_d:=\{(\gamma,\eta):|\eta^{-1}\gamma|_r=d\}\subset G(\Gamma)\times G(\Gamma)$ and let $\chi_d$ be the characteristic function for $E_d$. By summing $\psi_{n,\Gamma,d}$, over all $d\in \mathbb{N}$, we obtain a kernel $\psi_{n,\Gamma}$ on $G(\Gamma)$. 
\begin{equation*}
    \begin{array}{rl}
        \psi_{n,\Gamma}:G(\Gamma)\times G(\Gamma)  \rightarrow & \mathbb{C} \\
         (\gamma,\eta)   \mapsto & \langle Exp_o(R_\Gamma(\gamma)),Exp_o(R_\Gamma(\eta)) \rangle \sum_{d=0}^\infty e^{-\frac{d}{n}}\psi_{n,\Gamma,d}(\gamma,\eta)\chi_d(\gamma,\eta)\\
    \end{array}
\end{equation*}
In the following, we verify that these kernels $\psi_{n,\Gamma}$'s are completely bounded, and also prove them to be $G(\Gamma)$-invariant, inducing a sequence of functions $\phi_{n,\Gamma}$ on $G(\Gamma)$ with the same $B_2$-norm as that of the kernels.

\medskip

\subsection{$B_2$-norm estimate of $\psi_{n,\Gamma}$}\label{section II}

$R_\Gamma$ maps the graph product groups $G(\Gamma)$ into the Hilbert space $\mathcal{H}_1$ (see \ref{definition of R_Gamma}), and $(\gamma,\eta)\mapsto ||R_\Gamma(\gamma)-R_\Gamma(\eta)||^2$ is a CND kernel on $G(\Gamma)$. Therefore the kernel
\begin{equation*}
    (\gamma,\eta) \mapsto \langle Exp_o(\frac{R_\Gamma(\gamma)}{\sqrt{n}}),Exp_o(\frac{R_\Gamma(\eta)}{\sqrt{n}})\rangle=e^{-\frac{||R_\Gamma(\gamma)-R_\Gamma(\eta)||^2}{n}}
\end{equation*}
is PD (from lemma \ref{schoenberg}), taking value 1 on the diagonals, that implies the $B_2$-norm of the above kernel is 1.

\medskip

Recall that for any $x \in G_v$, the vector $avg_{n,v,S_v}(x)$ was defined to be an $\epsilon$-perturbed average of $\alpha_{n,v,S_v}(x)$ and $\beta_{n,v,S_v}(x)$. A simple calculation of norms will give us 

\begin{equation*}
  \left.\begin{aligned}
  ||\alpha_{n,v,S_v}(x)-avg_{n,v,S_v}(x)||^2&\leq\epsilon\\
  ||\beta_{n,v,S_v}(x)-avg_{n,v,S_v}(x)||^2&\leq\epsilon
\end{aligned}\right\}~ \forall ~ x \in G_v,
\end{equation*}
and
\begin{equation*}
    ||avg_{n,v,S_v}(x)||^2=\frac{||\alpha_{n,v,S_v}(x)||^2+||\beta_{n,v,S_v}(x)||^2+2}{(2+\sqrt{2\epsilon})^2}\leq \frac{2(2+\epsilon)}{(2+\sqrt{2\epsilon})^2},    
\end{equation*}
which in turn implies  $|D_{n,v,S_v}(x)|\geq\frac{\epsilon^{\frac{1}{4}}}{4}$. The above inequalities put together, provides a bound on the constants defined above as $C^{\alpha_{n,v,S_v}}(x,y)$ and $C^{\beta_{n,v,S_v}}(x,y)$, i.e. there exists a positive real number $A>0$ such that 

\begin{equation*}
    \sup_{x,y \in G_v} \{|C^{\alpha_{n,v,S_v}}(x,y)|,|C^{\beta_{n,v,S_v}}(x,y)|\} < A\epsilon^{\frac{1}{4}}
\end{equation*}

\medskip

In order to calculate the $B_2$-norm of $\psi_{n,\Gamma}$, we need to know what is the value of that norm for the kernel $\psi_{n,\Gamma,d}$. The above estimates gives us that the norm of the each component vectors of the vector $\alpha_{n,\Gamma,d}$ and $\beta_{n,\Gamma,d}$ respectively, is bounded by $(1+\epsilon+2A^2\epsilon^{\frac{1}{2}})^{\frac{1}{2}}<(1+B\sqrt{\epsilon})^{\frac{1}{2}}$, for some $B>0$.
Lemma \ref{maximum terms} tells that in the $d$-tail of any element $\gamma\in G(\Gamma)$, there can be at most $dM$ many letters. Hence the norms of the vectors $\alpha_{n,\Gamma,d}(\gamma)$ and $\beta_{n,\Gamma,d}(\gamma)$ are bounded by $(1+B\sqrt{\epsilon})^{\frac{dM}{2}}$ for any $\gamma \in G(\Gamma)$. Hence, from lemma \ref{cb equivalent}, 
\begin{equation*}
    ||\psi_{n,\Gamma,d}||_{B_2} < (1+B\sqrt{\epsilon})^{dM}
\end{equation*}
Recall that, by definition for any $(\gamma,\eta)\in G(\Gamma)\times G(\Gamma)$, we have 
\begin{equation}\label{sum kernel}
    \psi_{n,\Gamma}(\gamma,\eta)= \langle Exp_o(R_\Gamma(\gamma)),Exp_o(R_\Gamma(\eta))\rangle \sum_{d=0}^\infty e^{-\frac{d}{n}}\psi_{n,\Gamma,d}(\gamma,\eta)\chi_d(\gamma,\eta)
\end{equation}

\medskip

We know that the $B_2$-norm of the first part of the above kernel is 1. We will show that the summation in the second part is also close (in $B_2$-norm) to a PD kernel, whose $B_2$-norm is 1. This justifies the choices for the vectors $\alpha_{n,\Gamma,d}$ and $\beta_{n,\Gamma,d}$. 

\medskip

Suppose $\gamma=\gamma_1\gamma_2\cdots\gamma_m$ is a reduced form of $\gamma \in G(\Gamma)$. Define a map $\zeta_{n,\Gamma}:G(\Gamma)\rightarrow\mathcal{H}_2$ as follows:
$$\zeta_{n,\Gamma}(\gamma):=(\theta_{n,v_1,S_{v_1}}(\gamma_1))_{G(st(v_1))}\otimes\cdots \otimes(\theta_{n,v_m,S_{v_m}}(\gamma_m))_{\gamma_1\cdots\gamma_{m-1}G(st(v_m))}.$$
The well definedness of $\zeta_{n,\Gamma}$ is clear from the discussion of well definedness of the functions $\alpha_{n,\Gamma,d}$ given above and $||\zeta_{n,\Gamma}(\gamma)||=1$ for any element $\gamma\in G(\Gamma)$. Let us define the kernel $\sigma_{n,\Gamma}:G(\Gamma)\times G(\Gamma) \rightarrow \mathbb{C}$ as
\begin{equation}\label{positive deinite comparison}
    \sigma_{n,\Gamma}(\gamma,\eta)=\langle \zeta_{n,\Gamma}(\gamma),\zeta_{n,\Gamma}(\eta)\rangle
\end{equation}
and observe that it is PD with $B_2$-norm 1. We will take the help of this kernel $\sigma_{n,\Gamma}$ to estimate the $B_2$-norm of the second part of the equation \ref{sum kernel}. Notice that 
\begin{align*}
        \psi_{n,\Gamma,d}(\gamma,\eta)-\sigma_{n,\Gamma}(\gamma,\eta)=\langle\alpha_{n,\Gamma,d}
        (\gamma)& -\zeta_{n,\Gamma}(\gamma), \beta_{n,\Gamma,d}(\eta)\rangle\\ 
        & +\langle\zeta_{n,\Gamma}(\gamma),\beta_{n,\Gamma,d}(\eta)-\zeta_{n,\Gamma}(\eta)\rangle    
\end{align*}
Therefore
\begin{multline}\label{1st B_2 estimate}
    ||\psi_{n,\Gamma,d}-\sigma_{n,\Gamma}||_{B_2}\leq  \sup_{\gamma,\eta\in G(\Gamma)}\Big(||\alpha_{n,\Gamma,d}
        (\gamma) -\zeta_{n,\Gamma}(\gamma)||~||\beta_{n,\Gamma,d}(\eta)||+ \\
          ||\zeta_{n,\Gamma}(\gamma)||~||\beta_{n,\Gamma,d}(\eta)-\zeta_{n,\Gamma}(\eta)||\Big)
\end{multline}
Since for any $v\in \Gamma$ and $x \in G_v$, the following equations are true (from remark \ref{kernel at x,x}):
\begin{align*}
        \langle \alpha_{n,v,S_v}(x),avg_{n,v,S_v}(x) \rangle+ C^{\alpha_{n,v,S_v}}(x,x) \cdot  D_{n,v,S_v}(x)&=1 ~ ~ ~ ~\text{and} \\
    \langle \beta_{n,v,S_v}(x),avg_{n,v,S_v}(x) \rangle+ C^{\beta_{n,v,S_v}}(x,x) \cdot  D_{n,v,S_v}(x)&=1
\end{align*}
so that $\langle\alpha_{n,\Gamma,d}
        (\gamma),\zeta_{n,\Gamma}
        (\gamma)\rangle=1=\langle\zeta_{n,\Gamma}
        (\gamma),\alpha_{n,\Gamma,d}
        (\gamma)\rangle$ (same is true for $\beta_{n,\Gamma,d}$). Using this we have:
\begin{equation*}
    \begin{array}{rl}
         ||\alpha_{n,\Gamma,d}
        (\gamma) -\zeta_{n,\Gamma}(\gamma)||^2= & ||\alpha_{n,\Gamma,d}
        (\gamma)||^2-1\qquad \text{and} \\
        \\
         ||\beta_{n,\Gamma,d}(\gamma)-\zeta_{n,\Gamma}(\gamma)||^2= & ||\beta_{n,\Gamma,d}
        (\gamma)||^2-1 
    \end{array}
\end{equation*}
Hence the inequality \ref{1st B_2 estimate} becomes,
\begin{equation*}\label{B_2 norm compare}
    \begin{array}{rl}
        ||\psi_{n,\Gamma,d}-\sigma_{n,\Gamma}||_{B_2}\leq & 
        \sup_{\gamma,\eta\in G(\Gamma)}\Big(\sqrt{||\alpha_{n,\Gamma,d}(\gamma)||^2-1}~||\beta_{n,\Gamma,d}(\eta)|| \\
        
      & \qquad\qquad\qquad\qquad\qquad\qquad\qquad+~~\sqrt{||\beta_{n,\Gamma,d}(\eta)||^2-1} \Big) \\
        \\
        \leq & \sqrt{(1+B\sqrt{\epsilon})^{dM}-1} ~ (1+B\sqrt{\epsilon})^{\frac{dM}{2}}+ \sqrt{(1+B\sqrt{\epsilon})^{dM}-1}\\
        \\
        =& \Big(\sqrt{(1+B\sqrt{\epsilon})^{dM}-1}\Big)\cdot\Big((1+B\sqrt{\epsilon})^{\frac{dM}{2}}+1\Big)\\
        \\
        \leq & 2\sqrt{BdM}(1+B\sqrt{\epsilon})^{dM}\epsilon^{\frac{1}{4}}
    \end{array}
\end{equation*}

We now estimate $||\psi_{n,\Gamma}||_{B_2}$:
\begin{align*}
    ||\psi_{n,\Gamma}||_{B_2} &= ||\langle Exp_o(R_\Gamma(\cdot)),Exp_o(R_\Gamma(\cdot))\rangle \sum_{d=0}^\infty e^{-\frac{d}{n}}\psi_{n,\Gamma,d}(\cdot,\cdot)\chi_d(\cdot,\cdot)||_{B_2}\\
    &= ||\langle Exp_o(R_\Gamma(\cdot)),Exp_o(R_\Gamma(\cdot))\rangle||_{B_2} ~~||\sum_{d=0}^\infty e^{-\frac{d}{n}}\psi_{n,\Gamma,d}(\cdot,\cdot)\chi_d(\cdot,\cdot)||_{B_2}\\
    &=1\cdot ||\sum_{d=0}^\infty e^{-\frac{d}{n}}\psi_{n,\Gamma,d}(\cdot,\cdot)\chi_d(\cdot,\cdot)||_{B_2}
\end{align*}

\medskip

We need to estimate $B_2$-norm of $\chi_d$, since it appears in the above expression of $||\psi_{n,\Gamma}||_{B_2}$. Precisely this estimate was due to Mizuta in \cite{mizuta2008bozejko}. He showed the following:
\begin{lemma}\cite{mizuta2008bozejko}
    Suppose $X$ is a finite dimensional CAT(0) cube complex, and $\Xi_d$ is characteristic function of the set $\{(x,y):d_X(x,y)=d\}$ then there exists a polynomial $p$ such that $||\Xi_d||_{B_2}\leq p(d)$.
\end{lemma}

As mentioned in the last part of the section \ref{construction of wall space}, we have a finite dimensional CAT(0) cube complex $X$ on which $G(\Gamma)$ acts. The equation \ref{CAT} and the previous lemma together implies that $||\chi_d||_{B_2}$ on $G(\Gamma)$ is bounded by $p(d)$. So we can write,
\begin{align*}
    ||\psi_{n,\Gamma}||_{B_2} &\leq  ||\sum_{d=0}^\infty e^{-\frac{d}{n}}\Big(\psi_{n,\Gamma,d}(\cdot,\cdot)-\sigma_{n,\Gamma}(\cdot,\cdot)\Big)\chi_{d}(\cdot,\cdot)||_{B_2}+||\sum_{d=0}^\infty e^{-\frac{d}{n}}\sigma_{n,\Gamma}(\cdot,\cdot)\chi_{d}(\cdot,\cdot)||_{B_2}\\  
         &\leq \sum_{d=0}^\infty e^{-\frac{d}{n}}~||\psi_{n,\Gamma,d}-\sigma_{n,\Gamma}||_{B_2}~||\chi_d||_{B_2}+||\sigma_{n,\Gamma}||_{B_2}~||\sum_{i=0}^\infty e^{-\frac{d}{n}}\chi_d||_{B_2}\\
         &\leq ||\sigma_{n,\Gamma}||_{B_2}||e^{-\frac{|\cdot|_r}{n}}||_{B_2}+ \sum_{d=0}^\infty e^{-\frac{d}{n}}2\sqrt{BdM}(1+B\sqrt{\epsilon})^{dM}\epsilon^{\frac{1}{4}}p(d)
\end{align*}

\medskip

Here the first summand is 1, because $\sigma_{n,\Gamma}$ is PD with value 1 on diagonals, and the other term is also PD following equation \ref{CAT} and lemma \ref{schoenberg}. The second summand is bounded by $\epsilon^{\frac{1}{4}}\cdot\bigg(\sum_{d=0}^\infty e^{-\frac{d}{n}}(1+B\sqrt{\epsilon})^{dM}q(d)\bigg)$, where $q(d)$ is a polynomial of degree one higher than $p(d)$. Using ratio test, we can say that this series is convergent whenever $n\leq \frac{1}{M\ln{(1+B\sqrt{\epsilon})}}$.

\medskip

Let $\delta>0$ be a positive number, then for $n=\lfloor \frac{1}{M\ln{(1+B\sqrt{\epsilon})}} \rfloor$, there exists $K\in \mathbb{N}$, such that $\epsilon^{\frac{1}{4}}
\sum_{d=K+1}^\infty e^{-\frac{d}{n}}(1+B\sqrt{\epsilon})^{dM}q(d)<\frac{\delta}{2}$. Further we can choose $\epsilon$ sufficiently small so that the following is also true:

\begin{align*}
    \epsilon^{\frac{1}{4}}\bigg(\sum_{d=0}^K e^{-\frac{d}{n}}(1+B\sqrt{\epsilon})^{dM}q(d)\bigg) & \leq \epsilon^{\frac{1}{4}}(1+B\sqrt{\epsilon})^{KM}\sum_{d=0}^\infty e^{-\frac{d}{n}}q(d)<\frac{\delta}{2}.
\end{align*}

Then for that particular $n$ we will have $||\psi_{n,\Gamma}||_{B_2}<1+\delta$. Now for a decreasing sequence of $\epsilon_k\in (0,\epsilon)$, we have a subsequence of kernels $\{\psi_{n_k,\Gamma}\}$ with the property that $||\psi_{n_k,\Gamma}||_{B_2}<1+\delta$ for each $k$.




\subsection{The kernel $\psi_{n,\Gamma}$ is $G(\Gamma)$-invariant}\label{section III}
In this section we prove $G(\Gamma)$-invariance of $\psi_{n,\Gamma}$, i.e., for any three elements $\gamma,\eta,\eta' \in G(\Gamma)$ 
\begin{equation*}
    \psi_{n,\Gamma}(\eta,\eta')=\psi_{n,\Gamma}(\gamma\eta,\gamma\eta'),
\end{equation*}
which in turn defines the following function on $G(\Gamma)$:
\begin{align*}
    \phi_{n,\Gamma}:G(\Gamma) &\rightarrow \mathbb{C}\\
    \gamma & \mapsto \psi_{n,\Gamma}(\gamma,1_{G(\Gamma)}) 
\end{align*}
with the same $B_2$-norm as $\psi_{n,\Gamma}$. We do the following in order to get $G(\Gamma)$-invariance:
\begin{itemize}
    \item we show $\psi_{n,\Gamma}(\gamma,\eta)=\psi_{n,\Gamma}(\eta^{-1}\gamma,1_{G(\Gamma)})$, for any $\gamma,\eta \in G(\Gamma)$, and
    \item if $\psi_{n,\Gamma}(\gamma,\eta)=e^{f(\gamma,\eta)}$ for some function $f(\gamma,\eta)$ (see \ref{sum kernel}), for the $G(\Gamma)$-invariance of $\psi_{n,\Gamma}$, it is enough to show that $f(\gamma,\eta)=f(\eta^{-1}\gamma,1_{G(\Gamma)})$, for any $\gamma,\eta\in G(\Gamma)$.
\end{itemize}

\medskip

Let $\gamma,\eta$ be two arbitrary elements of $G(\Gamma)$, with two corresponding reduced forms $\gamma=\gamma_1\gamma_2\cdots\gamma_m$ and $\eta=\eta_1\eta_2\cdots\eta_p$. The reduced form configuration of $\eta^{-1}\gamma$ in $G(\Gamma)$, has the following possibilities:
\begin{enumerate}
    \item\label{1} the letters $ \gamma_i$ and $\eta_i $  are equal, i.e. $\gamma_i=\eta_i$ in $G_{v_i}$ for $1\leq i \leq q\leq \min\{m,p\}$
    \item\label{2} the letters $ \gamma_i$ and $\eta_i $ are in the same $G_{v_i}$ but are distinct for $1\leq i\leq r\leq \min\{m,p\}$, and the sub-graph in $\Gamma$ having vertices $\{v_{1},\cdots,v_{r}\}$ is a complete graph, where $r \leq M$ and $M$ is from the hypothesis of the lemma \ref{maximum terms},  
    \item\label{3} the possibilities described in (\ref{1}) and (\ref{2}) above may repeat finitely many times. These  possibilities, where $\gamma_i, \eta_i$ come from the same vertex group $G_{v_i}$, result in either cancellation or finite products of words of the form $\eta_i^{-1}\gamma_i$,
    \item\label{4}  there could be letters in the tail of the reduced forms of $\eta$ or $\gamma$ which survive in every reduced form of $\eta^{-1}\gamma$. 
\end{enumerate}
\begin{figure}[ht!]
    \centering
    \begin{tikzpicture}[scale= 0.45]
                \draw (-10,0) -- (-7,0) -- (-6,1) -- (-3,1) -- (-2,0) -- (0,0);
                \draw[dotted] (0,0) -- (2,0);
                \draw (2,0) -- (4,0) -- (5,1) -- (8,1) -- (9,2);
                \draw[dotted] (9,2) -- (11,4);
                \draw  (11,4)--(12,5);
                \draw (-7,0) -- (-6,-1) -- (-3,-1) -- (-2,0);
                \draw (4,0) -- (5,-1) -- (8,-1) -- (9,-2);
                \draw[dotted] (9,-2) -- (11,-4);
                \draw  (11,-4)--(12,-5);
                \draw [fill] (0,0) circle [radius=3pt];
                \draw [fill] (-10,0) circle [radius=3pt];
                \draw [fill] (-7,0) circle [radius=3pt];
                \draw [fill] (-6,1) circle [radius=3pt];
                \draw [fill] (-3,1) circle [radius=3pt];
                \draw [fill] (-2,0) circle [radius=3pt];
                \draw [fill] (2,0) circle [radius=3pt];
                \draw [fill] (4,0) circle [radius=3pt];
                \draw [fill] (5,1) circle [radius=3pt];
                \draw [fill] (8,1) circle [radius=3pt];
                \draw [fill] (-7,0) circle [radius=3pt];
                \draw [fill] (-6,-1) circle [radius=3pt];
                \draw [fill] (-3,-1) circle [radius=3pt];
                \draw [fill] (5,-1) circle [radius=3pt];
                \draw [fill] (8,-1) circle [radius=3pt];
                \draw [fill] (-9,0) circle [radius=3pt];
                \draw [fill] (-8,0) circle [radius=3pt];
                \draw [fill] (-5,1) circle [radius=3pt];
                \draw [fill] (-4,1) circle [radius=3pt];
                \draw [fill] (-5,-1) circle [radius=3pt];
                \draw [fill] (-4,-1) circle [radius=3pt];
                \draw [fill] (-1,0) circle [radius=3pt];
                \draw [fill] (3,0) circle [radius=3pt];
                \draw [fill] (6,1) circle [radius=3pt];
                \draw [fill] (6,-1) circle [radius=3pt];
                \draw [fill] (7,1) circle [radius=3pt];
                \draw [fill] (7,-1) circle [radius=3pt];
                \draw [fill] (9,2) circle [radius=3pt];
                \draw [fill] (11,4) circle [radius=3pt];
                \draw [fill] (12,5) circle [radius=3pt];
                \draw [fill] (9,-2) circle [radius=3pt];
                \draw [fill] (11,-4) circle [radius=3pt];
                \draw [fill] (12,-5) circle [radius=3pt];
                \path [draw=none,fill=gray, fill opacity = 0.2] (-6,0) ellipse (0.35cm and 1.5cm);
                \path [draw=none,fill=gray, fill opacity = 0.2] (-5,0) ellipse (0.35cm and 1.5cm);
                \node at (-5,-2.3) {$\Big\uparrow$};
                \node at (-5,-3.5) {some $G_v$-coset};
                \path [draw=none,fill=gray, fill opacity = 0.2] (-4,0) ellipse (0.35cm and 1.5cm);
                \path [draw=none,fill=gray, fill opacity = 0.2] (-3,0) ellipse (0.35cm and 1.5cm);
                \path [draw=none,fill=gray, fill opacity = 0.2] (5,0) ellipse (0.35cm and 1.5cm);
                \path [draw=none,fill=gray, fill opacity = 0.2] (6,0) ellipse (0.35cm and 1.5cm);
                \path [draw=none,fill=gray, fill opacity = 0.2] (7,0) ellipse (0.35cm and 1.5cm);
                \path [draw=none,fill=gray, fill opacity = 0.2] (8,0) ellipse (0.35cm and 1.5cm);
                \node at (-11.5,0) {$1_{G(\Gamma)}$};
                \node at (12.5,5.5) {$\eta$};
                \node at (12.5,-5.5) {$\gamma$};
                \node at (-8.5,-0.75) {$\underbrace{\qquad\qquad}$};
                \node at (-8.5,-1.5) {$\gamma_i=\eta_i$};
                \node at (-4.5,1.75) {$\overbrace{~~\qquad\qquad}$};
                \node at (-4.5,2.5) {$\eta_i^{-1}\gamma_i$};
                \node at (-1,-0.75) {$\underbrace{~~~~~\qquad}$};
                \node at (-1,-1.5) {$\gamma_i=\eta_i$};
                \node at (3,-0.75) {$\underbrace{~~~~~\qquad}$};
                \node at (3,-1.5) {$\gamma_i=\eta_i$};
                \node at (6.5,1.75) {$\overbrace{~~\qquad\qquad}$};
                \node at (6.5,2.5) {$\eta_i^{-1}\gamma_i$};
    \end{tikzpicture}
    \caption{Schematic diagram of reduced form of $\eta^{-1}\gamma$}
    \label{fig:4}
\end{figure}
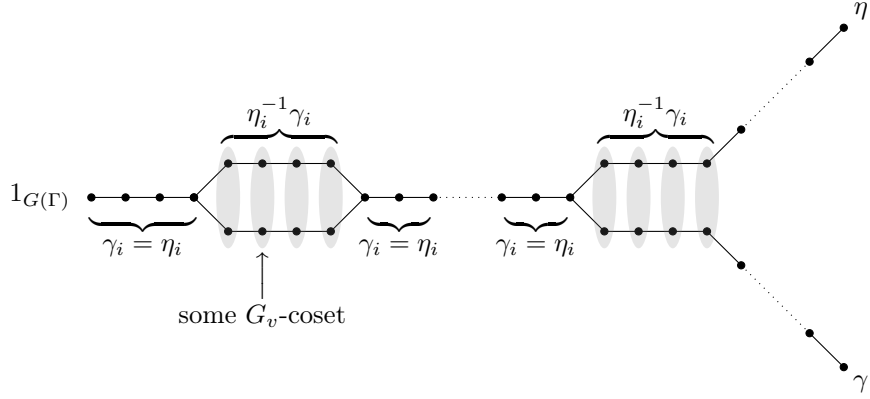

It is clear from the above discussion (Fig. \ref{fig:4}), that the reduced form of $\eta^{-1}\gamma$ consists of syllables or letters of types: (i) finite products of words of the form $\eta_i^{-1}\gamma_i$ (as the other possibilities mentioned in (\ref{1}) get cancelled), and (ii) the remaining letters surviving as mentioned in possibility \ref{4}. So the proof of $f(\gamma,\eta)=f(\eta^{-1}\gamma,1_{G(\Gamma)})$ gets reduce to the following case (see figure \ref{fig:5}):

Let $\gamma,\eta\in G(\Gamma)$ with a above mentioned reduced form such that $\gamma_i=\eta_i$ for $1\leq i \leq q\leq \min\{m,p\} $ and the sub-graph in $\Gamma$ having vertices $\{v_{q+1},\cdots,v_{q+r}\}$ is a complete graph with $r\leq M$, then a reduced form of $\eta^{-1}\gamma$ is 
\begin{equation*}
    \eta_p^{-1}\cdots\eta_{q+1+r}^{-1} (\eta_{q+1}^{-1}\gamma_{q+1})\cdots (\eta_{q+r}^{-1}\gamma_{q+r})   \gamma_{q+1+r}\cdots\gamma_m.
\end{equation*}
 We are going to calculate $f(\eta^{-1}\gamma,1_{G(\Gamma)})$ and $f(\gamma,\eta)$ separately, for this type of $\gamma$ and $\eta$.

\medskip

\subsubsection{Expression of $f(\eta^{-1}\gamma,1_{G(\Gamma)})$}

First we consider $\psi_{n,\gamma}(\eta^{-1}\gamma,1_{G(\Gamma)})$ and find out the exponent of $e$ in this expression:
\begin{align*}
    \psi_{n,\Gamma}(\eta^{-1}\gamma,1_{G(\Gamma)})=& \langle Exp_o(R_\Gamma(\eta^{-1}\gamma)),Exp_o(R_\Gamma(1_{G(\Gamma)}))\rangle\\ & \qquad\qquad\qquad \sum_{d=0}^\infty e^{-\frac{d}{n}}\psi_{n,\Gamma,d}(\eta^{-1}\gamma,1_{G(\Gamma)})\chi_d(\eta^{-1}\gamma,1_{G(\Gamma)})\\
    =&e^{-\frac{||R_\Gamma(\eta^{-1}\gamma)||^2}{n}}e^{-\frac{|\eta^{-1}\gamma|_r}{n}}\psi_{n,|\eta^{-1}\gamma|_r,\Gamma}(\eta^{-1}\gamma,1_{G(\Gamma)})
\end{align*}
Simplifying the above equation we get, 
\begin{multline*}
    \psi_{n,\Gamma}(\eta^{-1}\gamma,1_{G(\Gamma)})= e^{-\frac{||R_\Gamma(\eta^{-1}\gamma)||^2}{n}}e^{-\frac{|\eta^{-1}\gamma|_r}{n}}~\Pi_{i=q+1+r}^p\langle\alpha_{n,v_i,S_{v_i}}(\eta_i^{-1}),\beta_{n,v_i,S_{v_i}}(1_{v_i})\rangle\\
     \Pi_{i=q+1}^{q+r}\langle\alpha_{n,v_i,S_{v_i}}(\eta_i^{-1}\gamma_i),\beta_{n,v_i,S_{v_i}}(1_{v_i})\rangle\\
    \Pi_{i=q+1+r}^m\langle\alpha_{n,v_i,S_{v_i}}(\gamma_i),\beta_{n,v_i,S_{v_i}}(1_{v_i})\rangle
\end{multline*}
By substituting the values of the corresponding inner-products we get,
\begin{multline*}
    \psi_{n,\Gamma}(\eta^{-1}\gamma,1_{G(\Gamma)})= e^{-\frac{||R_\Gamma(\eta^{-1}\gamma)||^2}{n}}e^{-\frac{|\eta^{-1}\gamma|_r}{n}} e^{\sum_{i=q+1+r}^p\frac{\phi_{v_i}(1_{v_i})-||S_{v_i}(\eta_i^{-1})+S_{v_i}(1_{v_i})||^2}{n}}\\
    e^{\sum_{i=q+1}^{q+r}\frac{\phi_{v_i}(1_{v_i})-||S_{v_i}(\eta_i^{-1}\gamma_i)+S_{v_i}(1_{v_i})||^2}{n}} ~e^{\sum_{i=q+1+r}^{m}\frac{\phi_{v_i}(1_{v_i})-||S_{v_i}(\gamma_i)+S_{v_i}(1_{v_i})||^2}{n}}.
\end{multline*}
Hence, the total index of $e$, obtained from the right hand side of the above equation, is
\begin{multline*}
    -\frac{||R_\Gamma(\eta^{-1}\gamma)||^2}{n}-\frac{|\eta^{-1}\gamma|_r}{n} +\sum_{i=q+1+r}^p\frac{\phi_{v_i}(1_{v_i})-||S_{v_i}(\eta_i^{-1})+S_{v_i}(1_{v_i})||^2}{n}+\\ \qquad
    \sum_{i=q+1}^{q+r}\frac{\phi_{v_i}(1_{v_i})-||S_{v_i}(\eta_i^{-1}\gamma_i)+S_{v_i}(1_{v_i})||^2}{n} +
    \sum_{i=q+1+r}^{m}\frac{\phi_{v_i}(1_{v_i})-||S_{v_i}(\gamma_i)+S_{v_i}(1_{v_i})||^2}{n}
\end{multline*}
Therefore splitting $R_\Gamma$, by using its definition, the above expression becomes:
\begin{multline}\label{index calculation 1}
   f(\eta^{-1}\gamma,1_{G(\Gamma)}):= -\frac{|\eta^{-1}\gamma|_r}{n} +\sum_{i=q+1+r}^p\frac{\phi_{v_i}(1_{v_i})-||R_{v_i}(\eta_i^{-1})||^2-||S_{v_i}(\eta_i^{-1})+S_{v_i}(1_{v_i})||^2}{n} \\ +
    \sum_{i=q+1}^{q+r}\frac{\phi_{v_i}(1_{v_i})-||R_{v_i}(\eta_i^{-1}\gamma_i)||^2-||S_{v_i}(\eta_i^{-1}\gamma_i)+S_{v_i}(1_{v_i})||^2}{n}\\ +
    \sum_{i=q+1+r}^{m}\frac{\phi_{v_i}(1_{v_i})-||R_{v_i}(\gamma_i)||^2-||S_{v_i}(\gamma_i)+S_{v_i}(1_{v_i})||^2}{n}
\end{multline}

So, \ref{index calculation 1} gives the exponent on the R.H.S. of the equation: $$\psi_{n,\Gamma}(\eta^{-1}\gamma,1_{G(\Gamma)})=e^{f(\eta^{-1}\gamma,1_{G(\Gamma)})}.$$ 

\subsubsection{Expression of $f(\gamma,\eta)$}

Now we do the same thing in order to find out the exponent of $e$ in the value of $\psi_{n,\Gamma}(\gamma,\eta)$.
\begin{align*}
    \psi_{n,\Gamma}(\gamma,\eta)=& \langle Exp_o(R_\Gamma(\gamma)),Exp_o(R_\Gamma(\eta))\rangle \sum_{d=0}^\infty e^{-\frac{d}{n}}\psi_{n,\Gamma,d}(\gamma,\eta)\chi_d(\gamma,\eta)\\
    =&e^{-\frac{||R_\Gamma(\gamma)-R_\Gamma(\eta)||^2}{n}}e^{-\frac{|\eta^{-1}\gamma|_r}{n}}\psi_{n,|\eta^{-1}\gamma|_r,\Gamma}(\gamma,\eta)
\end{align*}
From the definition of $\psi_{n,|\eta^{-1}\gamma|_r,\Gamma}$ we have $$\psi_{n,|\eta^{-1}\gamma|_r,\Gamma}(\gamma,\eta)=\langle\alpha_{n,|\eta^{-1}\gamma|_r,\Gamma}(\gamma),\beta_{n,|\eta^{-1}\gamma|_r,\Gamma}(\eta)\rangle.$$

\medskip

Depending on how large the reduced length $|\eta^{-1}\gamma|_r$ is, the first $q$ terms of the inner product  \\$\langle\alpha_{n,|\eta^{-1}\gamma|_r,\Gamma}(\gamma),\beta_{n,|\eta^{-1}\gamma|_r,\Gamma}(\eta)\rangle$ become one of the following types (see the figure \ref{fig:5}):

\begin{figure}[ht!]
    \centering
        \begin{tikzpicture}[scale=0.3]
             \draw (-15,0) -- (-7,0) -- (-6,1) -- (1,1) -- (3,3);
             \draw[dotted] (3,3) -- (5,5);
             \draw  (5,5) -- (6,6);
             \draw (-7,0) -- (-6,-1) -- (1,-1) -- (3,-3);
             \draw[dotted] (3,-3) -- (5,-5);
             \draw  (5,-5) -- (6,-6);
             \draw [fill] (-15,0) circle [radius=3pt];
             \draw [fill] (-14,0) circle [radius=3pt];
             \draw [fill] (-13,0) circle [radius=3pt];
             \draw [fill] (-12,0) circle [radius=3pt];
             \draw [fill] (-11,0) circle [radius=3pt];
             \draw [fill] (-10,0) circle [radius=3pt];
             \draw [fill] (-9,0) circle [radius=3pt];
             \draw [fill] (-8,0) circle [radius=3pt];
             \draw [fill] (-7,0) circle [radius=3pt];
             \draw [fill] (-6,1) circle [radius=3pt];
             \draw [fill] (-5,1) circle [radius=3pt];
             \draw [fill] (-4,1) circle [radius=3pt];
             \draw [fill] (-3,1) circle [radius=3pt];
             \draw [fill] (-2,1) circle [radius=3pt];
             \draw [fill] (-1,1) circle [radius=3pt];
             \draw [fill] (0,1) circle [radius=3pt];
             \draw [fill] (1,1) circle [radius=3pt];
             \draw [fill] (-6,-1) circle [radius=3pt];
             \draw [fill] (-5,-1) circle [radius=3pt];
             \draw [fill] (-4,-1) circle [radius=3pt];
             \draw [fill] (-3,-1) circle [radius=3pt];
             \draw [fill] (-2,-1) circle [radius=3pt];
             \draw [fill] (-1,-1) circle [radius=3pt];
             \draw [fill] (0,-1) circle [radius=3pt];
             \draw [fill] (1,-1) circle [radius=3pt];
             \draw [fill] (2,2) circle [radius=3pt];
             \draw [fill] (2,-2) circle [radius=3pt];
             \draw [fill] (3,3) circle [radius=3pt];
             \draw [fill] (3,-3) circle [radius=3pt];
             \draw [fill] (5,5) circle [radius=3pt];
             \draw [fill] (6,6) circle [radius=3pt];
             \draw [fill] (5,-5) circle [radius=3pt];
             \draw [fill] (6,-6) circle [radius=3pt];
             \path [draw=none,fill=gray, fill opacity = 0.2] (-6,0) ellipse (0.35cm and 1.5cm);
             \path [draw=none,fill=gray, fill opacity = 0.2] (-5,0) ellipse (0.35cm and 1.5cm); 
             \path [draw=none,fill=gray, fill opacity = 0.2] (-4,0) ellipse (0.35cm and 1.5cm);
             \path [draw=none,fill=gray, fill opacity = 0.2] (-3,0) ellipse (0.35cm and 1.5cm); 
             \path [draw=none,fill=gray, fill opacity = 0.2] (-2,0) ellipse (0.35cm and 1.5cm); 
             \path [draw=none,fill=gray, fill opacity = 0.2] (-1,0) ellipse (0.35cm and 1.5cm);
             \path [draw=none,fill=gray, fill opacity = 0.2] (0,0) ellipse (0.35cm and 1.5cm); 
             \path [draw=none,fill=gray, fill opacity = 0.2] (1,0) ellipse (0.35cm and 1.5cm);
             \node at (-11,-0.75) {$\underbrace{\quad\qquad\qquad\qquad}$};
             \node at (-11,-2) {\footnotesize $\gamma_i=\eta_i$};
             \node at (-11,-3) {\footnotesize{$1\leq i \leq q$}};
             \node at (-2.5,1.75) {$\overbrace{\quad\quad\qquad\qquad}$};
             \node at (-2.5,4.5) {\footnotesize$\eta_i^{-1}\gamma_i$};
             \node at (-2.5,3) {\footnotesize{$q+1\leq i \leq q+r$}};
             \node at (6.5,-6.5) {$\gamma$};
             \node at (6.5,6.5) {$\eta$};
             \node at (-16.75,0) {$1_{G(\Gamma)}$};
        \end{tikzpicture}
    \caption{}
    \label{fig:5}
\end{figure}

\begin{align}\label{error manipulation}
    \begin{cases}
    1. & \langle \theta_{n,v_i,S_{v_i}}(\gamma_i),\theta_{n,v_i,S_{v_i}}(\eta_i)  \rangle\\
    2. & \langle \theta_{n,v_i,S_{v_i}}(\gamma_i),\Bigg(\beta_{n,v_i,S_{v_i}}(\eta_i),\begin{pmatrix}
        0\\C^{\beta_{n,v_i,S_{v_i}}}(\eta_i,\eta_i)
    \end{pmatrix},\begin{pmatrix}
        0\\C^{\beta_{n,v_i,S_{v_i}}}(\eta_i,1_{v_i})
    \end{pmatrix}\Bigg)\rangle\\
    3. & \langle \Bigg(\alpha_{n,v_i,S_{v_i}}(\gamma_i),\begin{pmatrix}
        C^{\alpha_{n,v_i,S_{v_i}}}(\gamma_i,\gamma_i)\\0
    \end{pmatrix},\begin{pmatrix}
        C^{\alpha_{n,v_i,S_{v_i}}}(\gamma_i,1_{v_i})\\0
    \end{pmatrix}\Bigg),\theta_{n,v_i,S_{v_i}}(\eta_i)\rangle\\
    4. &  \langle \Bigg(\alpha_{n,v_i,S_{v_i}}(\gamma_i),\begin{pmatrix}
        C^{\alpha_{n,v_i,S_{v_i}}}(\gamma_i,\gamma_i)\\0
    \end{pmatrix},\begin{pmatrix}
        C^{\alpha_{n,v_i,S_{v_i}}}(\gamma_i,1_{v_i})\\0
    \end{pmatrix}\Bigg),\\
    &\qquad\qquad\qquad\Bigg(\beta_{n,v_i,S_{v_i}}(\eta_i),\begin{pmatrix}
        0\\C^{\beta_{n,v_i,S_{v_i}}}(\eta_i,\eta_i)
    \end{pmatrix},\begin{pmatrix}
        0\\C^{\beta_{n,v_i,S_{v_i}}}(\eta_i,1_{v_i})
    \end{pmatrix}\Bigg)  \rangle.
    \end{cases}
\end{align}
In all these cases since $\gamma_i=\eta_i$, therefore the inner-product is 1. Lemma \ref{d tail lemma}, says that each $\gamma_i$ and $\eta_i$ are in the $|\eta^{-1}\gamma|_r$-tail of $\gamma$ and $\eta$ respectively, for $i\geq q+1$. Again from the given reduced form of $\eta^{-1}\gamma$, we can say that $\gamma_i,\eta_i \in G_{v_i}$ for $q+1\leq i\leq q+r$. Therefore the next $r$ terms are of the type (see the middle part of the figure \ref{fig:5}):
\begin{multline*}
    \langle \Bigg(\alpha_{n,v_i,S_{v_i}}(\gamma_i),\begin{pmatrix}
        C^{\alpha_{n,v_i,S_{v_i}}}(\gamma_i,\gamma_i)\\0
    \end{pmatrix},\begin{pmatrix}
        C^{\alpha_{n,v_i,S_{v_i}}}(\gamma_i,1_{v_i})\\0
    \end{pmatrix}\Bigg),\\
   \Bigg(\beta_{n,v_i,S_{v_i}}(\eta_i),\begin{pmatrix}
        0\\C^{\beta_{n,v_i,S_{v_i}}}(\eta_i,\eta_i)
    \end{pmatrix},\begin{pmatrix}
        0\\C^{\beta_{n,v_i,S_{v_i}}}(\eta_i,1_{v_i})
    \end{pmatrix}\Bigg)  \rangle
\end{multline*}
whose value is $\langle\alpha_{n,v_i,S_{v_i}}(\gamma_i),\beta_{n,v_i,S_{v_i}}(\eta_i)$. And, finally for $i\geq q+1+r$, the remaining terms of $\langle\alpha_{n,|\eta^{-1}\gamma|_r,\Gamma}(\gamma),\beta_{n,|\eta^{-1}\gamma|_r,\Gamma}(\eta)\rangle$ are of the form:

\begin{align*}
    \begin{cases}
        1. & \langle \Bigg(\alpha_{n,v_i,S_{v_i}}(\gamma_i),\begin{pmatrix}
        C^{\alpha_{n,v_i,S_{v_i}}}(\gamma_i,\gamma_i)\\0
    \end{pmatrix},\begin{pmatrix}
        C^{\alpha_{n,v_i,S_{v_i}}}(\gamma_i,1_{v_i})\\0
    \end{pmatrix}\Bigg),\theta_{n,v_i,S_{v_i}}(1_{v_i})\rangle\\
    &~ ~ ~=\langle\alpha_{n,v_i,S_{v_i}}(\gamma_i),\beta_{n,v_i,S_{v_i}}(1_{v_i})\rangle\\
    2. & \langle \theta_{n,v_i,S_{v_i}}(1_{v_i}),\Bigg(\beta_{n,v_i,S_{v_i}}(\eta_i),\begin{pmatrix}
        0\\C^{\beta_{n,v_i,S_{v_i}}}(\eta_i,\eta_i)
    \end{pmatrix},\begin{pmatrix}
        0\\C^{\beta_{n,v_i,S_{v_i}}}(\eta_i,1_{v_i})
    \end{pmatrix}\Bigg)\rangle\\
    & ~~~=\langle \alpha_{n,v_i,S_{v_i}}(1_{v_i}),\beta_{n,v_i,S_{v_i}}(\eta_i)\rangle
    \end{cases}
\end{align*}
So the value of $\psi_{n,\Gamma}(\gamma,\eta)$ is 
\begin{multline*}
    \psi_{n,\Gamma}(\gamma,\eta)= e^{-\frac{||R_\Gamma(\gamma)-R_\Gamma(\eta)||^2}{n}-\frac{|\eta^{-1}\gamma|_r}{n}}e^{\sum_{i=q+1}^{q+r} \frac{\phi_{v_i}(1_{v_i})-||S_{v_i}(\gamma_i)+S_{v_i}(\eta_i)||^2}{n} } \cdot \\
     e^{\sum_{i=q+1+r}^{m}\frac{\phi_{v_i}(1_{v_i})-||S_{v_i}(\gamma_i)+S_{v_i}(1_{v_i})||^2}{n}}\cdot e^{\sum_{i=q+1+r}^{p}\frac{\phi_{v_i}(1_{v_i})-||S_{v_i}(1_{v_i})+S_{v_i}(\eta_i)||^2}{n} }
\end{multline*}
Again, by considering the index of $e$ and simplifying the expression of $R_\Gamma$, we get:
\begin{multline}\label{index calculation 2}
   f(\gamma,\eta):= -\frac{|\eta^{-1}\gamma|_r}{n}+\sum_{i=q+1}^{q+r}\frac{\phi_{v_i}(1_{v_i})-||R_{v_i}(\gamma_i)-R_{v_i}(\eta_i)||^2-||S_{v_i}(\gamma_i)+S_{v_i}(\eta_i)||^2}{n} +\\
  \sum_{i=q+1+r}^{m}\frac{\phi_{v_i}(1_{v_i})-||R_{v_i}(\gamma_i)||^2-||S_{v_i}(\gamma_i)+S_{v_i}(1_{v_i})||^2}{n}+\\
  \sum_{i=q+1+r}^p\frac{\phi_{v_i}(1_{v_i})-||R_{v_i}(\eta_i)||^2-||S_{v_i}(1_{v_i})+S_{v_i}(\eta_i)||^2}{n}
\end{multline}
From equations \ref{index calculation 1} and \ref{index calculation 2} it is clear that $f(\eta^{-1}\gamma,1_{G(\Gamma)})=f(\gamma,\eta)$, proving the function $\psi_{n,\Gamma}$ is $G(\Gamma)$-invariant, i.e. for any $\gamma, \eta,\eta^{\prime}\in G(\Gamma)$:
\begin{equation*}
    \psi_{n,\Gamma}(\eta,\eta')=\psi_{n,\Gamma}(\gamma\eta,\gamma\eta'),
\end{equation*}

\subsection{Proof of the main theorem}

In short, what we have from the above sections are:
\begin{itemize}
    \item in section \ref{section I}, we constructed a sequence of kernels $\psi_{n,\Gamma}$,
    \item section \ref{section II} gives a bound on a subsequence $||\psi_{n_k,\Gamma}||_{B_2}$,
    \item section \ref{section III} proves the invariance of $\psi_{n,\Gamma}$ under the diagonal action of $G(\Gamma)$.
\end{itemize}

\medskip

Hence the function $\phi_{n,\Gamma}:G(\Gamma)\rightarrow\mathbb{C}$ given by $\phi_{n,\Gamma}(\gamma):=\psi_{n_n,\Gamma}(\gamma,1_{G(\Gamma)})$ has $B_2$-norm bounded by $1+\delta$ and for a reduced form $\gamma_1\cdots\gamma_m$ of $\gamma$, the map looks like: 
\begin{equation*}
    \phi_{n,\Gamma}(\gamma)=e^{-\frac{|\gamma|_r+\sum_{i=1}^m\Big(\phi_{v_i}(\gamma_i)-\phi_{v_i}(1_{v_i})\Big)}{n}}
\end{equation*}
The map $\gamma \mapsto |\gamma|_r+\sum_{i=1}^m\Big(\phi_{v_i}(\gamma_i)-\phi_{v_i}(1_{v_i})\Big)$ is a proper map on $G(\Gamma)$, (see lemma 4.4, \cite{das2023stability} for a proof of properness). This proves that $\phi_{n,\Gamma}$ is a function vanishing at infinity.

\medskip

Summarizing what we have obtained so far: for any given $\delta>0$, we have maps 
\begin{equation*}
    \phi_{n,\Gamma}:G(\Gamma) \rightarrow \mathbb{C}
\end{equation*}
having properties:
\begin{enumerate}
    \item $\phi_{n,\Gamma}$ vanishes at infinity,
    \item $\phi_{n,\Gamma}(\gamma)\rightarrow 1$ as $n \rightarrow \infty$, for all $\gamma\in G(\Gamma)$,
    \item $||\phi_{n,\Gamma}||_{B_2}\leq 1+\delta$
\end{enumerate}
Hence we have proved the following theorem.

\begin{theorem}\label{main theorem}
    Suppose $\Gamma$ is a finite graph and for each vertex $v \in V(\Gamma)$, the associated group $G_v$ has the weak Haagerup property with $\Lambda_{WH}(G_v)=1$. Then $G(\Gamma)$ has weak Haagerup property. Moreover $\Lambda_{WH}(G(\Gamma))=1$.
\end{theorem}

The following is a special case:

\begin{corollary}\label{main corollary}
    Suppose $A$ and $B$ are two weakly Haagerup groups such that $\Lambda_{WH}(A)=1=\Lambda_{WH}(B)$, then $G=A \ast B$ has weak Haagerup property with $\Lambda_{WH}(G)=1$.
\end{corollary}

\begin{remark}
        The proof of corollary \ref{main corollary} will be same as theorem \ref{main theorem}. But here we take $X$ to be the Bass-Serre tree, which is a special case of a CAT(0)-cube complex, associated with the free products of groups \cite{serre2002trees}. It was proven in \cite{haagerup1978example}, that the distance metric on a tree or more precisely the normal length on $G$ is CND type. Finally it is shown by Bozejko and Picardello, in \cite{bozejko1993weakly}, that $||\chi_d||_{B_2}$ (appearing in $\psi_{n,\Gamma}$) for free products has a linear upper-bound i.e. there exists some $D>0$ such that $||\chi_d||_{B_2} \leq D(d+1)$.
    \end{remark}

\section*{Acknowledgment}
PSG\footnote{parthasarathi.ghosh.100@gmail.com} was supported by CSIR Fellowship [File No. 08/155(0066)/2019-EMR-I], Govt. of India. SD\footnote{shubhabrata.maths@presiuniv.ac.in} acknowledges the infrastructural support provided by the Dept. of Mathematics, Presidency University, through the DST-FIST [File No. SR/FST/MS-I/2019/41]. 

\bibliographystyle{plain} 
\bibliography{references}

\end{document}